\documentclass[11pt, oneside]{amsart}
\usepackage[centertags]{amsmath}
\usepackage{multicol}
\usepackage{amsfonts}
\usepackage{amssymb}
\usepackage{amsthm}

\usepackage{footnpag}

\usepackage{amscd}
\usepackage{amssymb}

\usepackage{curves}

\def\osc{\text{osc}}

\def\vz{\varepsilon}

\def\dz{\delta}

\def\wt{\tilde}

\def\sub{\substack}

\def\bi{\bibitem}

\newcommand{\al}{\alpha}
\def\be{\beta}
\def\sph{\mathbb{S}^{d}}
\def \b{\beta}

\def\bl{\bigl}
\def\br{\bigr}

\def\dmin{\displaystyle \min}
\def\Ld{\Lambda}
\def\da{\delta}

\def\df{\displaystyle\frac}

\def\dsup{\displaystyle\sup}
\def\dsum{\displaystyle\sum}

\def\og{\omega}

\def\sub{\substack}
\def\Bl{\Bigl}
\def\Br{\Bigr}
\def\f{\frac}
\def\({\Bigl(}
\def \){ \Bigr)}

\def\ga{\gamma}

\def\hb{\hfill$\Box$}

\newcommand{\ld}{\lambda}

\newcommand\p{\partial}

\def\sa{\sigma}
\def\dmax{\displaystyle\max}

\def\ta{\theta}
\newcommand{\tr}{\triangle}

\def\va{\varepsilon}

\def\wt{\widetilde}

\theoremstyle{plain}
\newtheorem{thm}{Theorem}[section]
\newtheorem{cor}[thm]{Corollary}
\newtheorem{lem}[thm]{Lemma}

\theoremstyle{definition}

\theoremstyle{remark}

\numberwithin{equation}{section}

\newfam\msbfam
\font\tenmsb=msbm10    \textfont\msbfam=\tenmsb \font\sevenmsb=msbm7
\scriptfont\msbfam=\sevenmsb \font\fivemsb=msbm5
\scriptscriptfont\msbfam=\fivemsb
\begin{document}
\title[]{ Positive Cubature formulas and\\ Marcinkiewicz-Zygmund
inequalities \\on spherical caps}
\author{Feng Dai}
\address{Department of Mathematical and Statistical Sciences,
University of Alberta, Edmonton, Alberta, T6G 2G1,  Canada. }
\email{ dfeng@math.ualberta.ca.}

\author{ Heping Wang} \address{ Department of Mathematics, Capital Normal
University,
\\Beijing 100037,
 People's Republic of China.}
 \thanks{ Research started while the second author visited Edmonton.  The first author  was partially supported  by the
NSERC Canada under grant G121211001.  The second  author was
partially  supported
 by the NNSF  China
(No. 10201021), the  SRCP Beijing (No. KM200310028106), and the
 NSF Beijing (No.  1062004).}

\keywords{ Spherical caps,   cubature formulas,
Marcinkiewicz-Zygmund inequalities,  spherical polynomials. }

\begin{abstract}
Let  $\Pi_n^d$ denote  the space of all spherical polynomials of
degree at most $n$ on the unit sphere $\sph$ of
$\mathbb{R}^{d+1}$, and let $d(x, y)$ denote the usual geodesic
distance $\arccos x\cdot y$ between $x, y\in \sph$.   Given a
spherical cap
$$ B(e,\al)=\{ x\in\sph:\   \ d(x, e)  \leq \al\}, \     \     \    \  (
e\in\sph,\   \  \text{$\al\in (0,\pi)$ is bounded away from
$\pi$}),$$
 we define the metric  $$\rho(x,y):=\frac 1{ \al}
\sqrt{(d(x, y))^2+\al(\sqrt{\al-d(x, e)}-\sqrt{\al-d(y,e)})^2},
  $$
  where  $x, y\in B(e,\al)$.
  It is shown  that given any $\be\ge 1$,
$1\leq p<\infty$ and any finite  subset $\Ld$  of $B(e,\al)$
satisfying the condition $\dmin_{\sub{\xi,\eta \in \Ld\\\xi\neq
\eta}} \rho (\xi,\eta) \ge \f \da n$ with $\da\in (0,1]$, there
exists a positive constant $C$, independent of $\al$, $n$, $\Ld$
and $\da$,  such that,   for any $f\in\Pi_{n}^d$,
 \begin{equation*} \sum _{\og\in \Ld} \(\max_{x,y\in
B_\rho (\og,\, \be\da/n)}|f(x)-f(y)|^p\)\,|B_\rho(\og, \da/n)| \le
(C \dz)^p \int_{B(e,\al)} |f(x)|^p\,d\sa(x),\end{equation*} where
$d\sa(x)$ denotes the usual Lebesgue measure on $\sph$,
 $$B_\rho(x, r )=\Bl\{y\in B(e,\al):\ \   \rho(y,x)\leq
 r\Br\},\    \    \  (r>0),$$
  and
 $$\Bl|B_\rho(x, \f\da n )\Br|=\int_{B_{\rho}(x, \da/n)} d\sa(y)
 \sim \al ^{d}\Bl[ \(\f{\da }n \)^{d+1}+ \(\f\da n\)^{d} \sqrt{1-\f{d(x, e)}\al}\Br].$$ As a consequence, we establish  positive
cubature formulas  and Marcinkiewicz-Zygmund inequalities on the
spherical cap $B(e,\al)$.
\end{abstract}
\maketitle

\newpage

\section{Introduction}

Let $\sph=\{x=(x_1,x_2,\cdots, x_{d+1})\in\mathbb{R}^{d+1}:\ \
|x|:=\sqrt{x_1^2+x_2^2+\cdots+x_{d+1}^2}=1\}$ denote   the unit
sphere of $\mathbb{R}^{d+1}$ endowed with the usual rotation
invariant measure $d\sa(x)$.  We denote by $d(x,y)$ the geodesic
distance $\arccos x\cdot y$ between $x$ and $y$ on $\sph$,  by
$B(x,r)$ the spherical cap $\{ y\in\sph:\ \ d(x,y)\leq r\}$
centered at $x\in\sph$ of radius $r>0$, and by $B(x; \al, \al+\b)$
the spherical collar $\{y\in\sph:\ \ \al\leq d(x,y)\leq \al+\b\}$
centered at $x\in\sph$ of spherical height $\b>0$. A function on
$\sph$ is called  a spherical polynomial  of degree at most $n$ if
it is the restriction to $\sph$ of a polynomial in $d+1$ variables
of total degree at most $n$.  We denote by $\Pi_n^d$
 the space of all spherical polynomials of degree at most $n$ on
$\sph$.   Given a set $E$, we shall use the notations   $\# E$,
$|E|$ and $\chi_E$ to denote  its  cardinality, Lebesgue measure,
and    characteristic function, respectively. Moreover, we shall
write $A\sim B$ for the statement $C^{-1} \leq A/B \leq C$, where
$C>0$ is called the constant of equivalence.

Let   $e$ be a fixed point on
 $\sph$,  $\al \in (0, \pi)$, $1\leq p<\infty$, $n$ be a positive
integer, and let $\Ld$ be a finite subset of the spherical cap
$B(e,\al)$.  We are particularly interested in the following
question: \vspace{2mm}

{\bf Question:}\  \ {\it  What condition on the finite subset
$\Ld$ guarantees  the existence of   a sequence of positive
numbers $\ld_\og$, $\og\in \Ld$  for which the following two
equations hold  for all $ f\in \Pi^d_n$?
\begin{equation}\label{1-1}
\int_{B(e,\al)} f(y)\, d\sa(y) =\sum_{ \og \in \Ld} \ld_\og f(
\og)\end{equation} and \begin{equation}\label{1-2} \int_{B(e,\al)}
|f(x)|^p\, d\sa(x) \sim \sum_{\og\in\Ld} \ld_\og |f(\og)|^p
,\end{equation} where the constant of equivalence is independent
of $n$,  $f$, and $\al$.}

\vspace{2mm}

In  the above question, we would expect  a sharp estimate on the
weights $\ld_\og$ and that $\#\Ld \sim \text{dim}\, \Pi_n^d$ as
$n\to \infty$. It should be pointed out that our interest in the
above question is motivated by the work [M1] of Mhaskar and the
work [KL2, GLN] of Golinskii, Kobindarajah,   Lubinsky and  Nevai.

 Note that the integrals in ($\ref{1-1}$)  and
($\ref{1-2}$) are over the spherical cap $B(e,\al)$ rather than
the whole sphere $\sph$.
 An equality like ($\ref{1-1}$) with positive
weights  $\ld_\og$ is called a positive cubature formula of degree
$n$, while an equivalence like ($\ref{1-2}$) is called a
Marcinkiewicz-Zygmund (MZ) type inequality.

In one dimensional case, MZ inequalities over arcs of the circle
for the full range of $0<p<\infty$ were obtained  in the paper
[GLN, Theorem 1.1] of Golinskii, Lubinsky, and Nevai, and in a
more recent  paper [KL2] of  Kobindarajah and   Lubinsky  (see
Remark 1.5 below for more details).   In  the case of $d\ge 2$,
existence of positive cubature formulas based on scattered data on
spherical caps were proved  by  Mhaskar in  [M1] (see Remark 1.6
below for details.) Many useful cubature formulas on $\sph$ with
different properties were previously
 constructed  by many authors.
For relevant results on $\sph$,  one may consult  [BD, BDS, D2,
GS, MNW, M2, NPW1, NPW2,  Pe,  Xu1, Xu2], among others.

Here we   recall some known results on  $\sph$ that are relevant
to our current discussion.  Positive cubature formulas and MZ
inequalities for $1\leq p\leq \infty$ based on function values at
scattered sites on  $\sph$ were first established by Mhaskar,
Narcowich and  Ward in the fundamental paper [MNW]. Positive
cubature formulas  on $\sph$ with  sharp estimates  on the weights
were  obtained by
 Narcowich, Petrushev and  Ward in a more recent paper [NPW1]. A
different proof of  MZ inequalities   on $\sph$   was given in
[BD].    Compared with that of [MNW], the proof in [BD] works for
the full range of $0<p<\infty$ and  for all compact two-point
homogeneous manifolds. For results concerning   doubling weights
on $\sph$, we refer to [D1]  (in the case $d\ge 2$) and the
remarkable work [MT2] of  Mastroianni and  Totik (in the case
$d=1$).

   It should be pointed out that all known  proofs of the MZ inequalities on $\sph$  are   based on
   the following integral representation of spherical polynomials:
\begin{equation}\label{1-3-re}f(x)=\int_{\sph} f(y) K_n (x\cdot y)\, d\sa(y),\ \ x\in
\sph,\   \ f\in \Pi_n^d,\end{equation} where $K_n$ is a smooth
reproducing kernel for the space $\Pi_n^d$  ( see, for instance,
[BD, (2.13)]). Since  the integral in  (\ref{1-3-re})  is over the
whole sphere rather than  on a local spherical cap,  we find it
difficult to  use (\ref{1-3-re}) to deduce  similar results  on
local spherical caps.
 In our opinion,
in order to obtain an ideal  result on a spherical cap $B(e,\al)$,
special efforts have to be made to treat  the center $e$ as well
as  the boundary of $B(e,\al)$. Our proof will be different from
those for $\sph$ (see, for instance, [MNW, BD]). It is based on
some recent results obtained in [BD] and [D1], as well as  the
weighted Markov-Bernstein-type inequality recently proved by Erd\`
elyi [Er2], rather than the integral representation
(\ref{1-3-re}).

To  state our main results, we need to introduce  several
necessary notation. Let $(X, d_X)$ be a metric space. We denote by
$B_{_{d_X}}(x, r)$ the ball $\{ y\in X: \   \ d_X (x, y) \leq r\}$
centered at $x\in X$ of radius $r>0$. Given $\va>0$ and a finite
subset $A$ of $X$, we say $A$ is $(\va, d_X)$-separable if it
satisfies the condition
$ \dmin_{\sub{\xi, \xi'\in A\\
\xi\neq\xi'}} d_X (\xi,\xi')\ge \vz;$    while  we say  $A$ is
maximal $(\va, d_X)$-separable if it is  $(\va, d_X)$-separable
and satisfies $ \displaystyle X=\bigcup_{\xi\in
A}B_{{_{d_X}}}(\xi,\vz).$

For $x, y\in B(e,\al)$, we define
\begin{equation}\label{1-3}
\rho(x,y)\equiv \rho_{_{B(e,\al)}}(x,y):=\frac 1{ \al}
\sqrt{(d(x,y))^2+\al(b_x^{1/2}-b_y^{1/2})^2},\end{equation} where
 $b_x\equiv b_{x, B(e,\al)}$    denotes  the shortest
distance from $x\in B(e,\al)$ to the boundary of  $B(e, \al)$;
that is \begin{equation}\label{1-4} b_x\equiv b_{x,
B(e,\al)}=\al-d(x,e).\end{equation} It's easily seen that $\rho$
is a metric on $B(e,\al)$. For $r>0$ and $x\in B(e,\al)$, we
define
\begin{equation}\label{1-5} \tr_r(x)\equiv \tr_{r, B(e,\al)}(x) :=
\al ^{d}\( r^{d+1}+ r^{d} \sqrt{1-\f{d(x, e)}\al}\).\end{equation}
 It will be  shown  in
Section 2 (Lemma 2.2 (iii)) that for any $x\in B(e,\al)$ and $r\in
(0,1)$,
$$ |B_\rho (x, r)| \sim \tr_{r}(x),$$
where and throughout $B_\rho(x,r):=\{ y\in B(e,\al):\   \ \rho(y,
x)\leq r\}$, and  the constant of equivalence is independent of
$r$, $x$, and $\al$ when $\al$ is bounded away from  $\pi$.

For the rest of this section,  we assume  that  $B(e, \al)$ is
given with  $\al\in (0,\pi)$ bounded away from $\pi$,  and we
write $\rho$, $b_x$ and $\tr_r(x)$ for $\rho_{_{B(e,\al)}}$,
$b_{x, B(e,\al)}$ and $\tr_{r, B(e,\al)}(x)$ respectively.

Now our main  result in this paper   can  be stated as follows.

\begin{thm}
If   $\da\in (0,1]$, $\be \ge 1$,   $1\leq p<\infty$, and $\Ld$ is
a $(\da/n, \rho)$-separable subset of $B(e, \al)$, then for any
$f\in\Pi_{n}^d$,  we have
\begin{equation}\label{1-6} \sum _{\og\in \Ld} \(\max_{x,y\in
B_\rho (\og,\, \be\da/n)}|f(x)-f(y)|^p\)\,|B_\rho(\og, \da/n)| \le
(C_1 \dz)^p \int_{B(e,\al)} |f(x)|^p\,d\sa(x),\end{equation} where
$C_1$ depends only on $d$, $p$, and  $\be$.
\end{thm}

 It can be shown that any maximal $(\da/n, \rho)$-separable
subset $\Ld$ of $B(e,\al)$ must satisfy the condition  $ \#\Ld
\sim \da^{-d} \text{dim} \, \Pi_n^d \sim \( \f n\da\)^d $, with
the constants of equivalence depending only on $d$. In particular,
this means that the number of nodes required in the above theorem
is comparable to the dimension of $\Pi_n^d$ as $n\to \infty$.

As a consequence of  Theorem 1.1, we have   the following three
useful corollaries.

\begin{cor} There exists a constant $\da_0\in (0,1)$ depending only on $d$,  such that,  for  any $\da \in (0, \da_0)$ and any maximal $(\da/n,
\rho)$-separable subset $\Ld$ of $B(e,\al)$, there exists a
sequence of positive numbers $\ld_\og$, $\og \in\Ld$ satisfying
\begin{equation}\label{1-7-a}\ld_\og  \sim \tr_{_{\da/n}} (\og)\sim \Bl|B_\rho
(\og, \f \da n)\Br|, \ \   \og\in\Ld\end{equation} with constants
of equivalence  depending  only on $d$, such that the  cubature
formula ($\ref{1-1}$) holds for all  $f\in \Pi_n^d$.
\end{cor}

Corollary 1.2 seems  new even in the case $d=1$. It can be deduced
from Theorem 1.1 following  the standard method  in [NPW1].

\begin{cor} Given  $1\leq p <\infty$ and an arbitrary  finite subset $\Ld$  of $B(e,\al)$,  there exists
a positive  constant  $C$ depending only on $p$ and $d$, such
that, for any  $f\in\Pi_n^d$,
$$ \sum_{\og \in \Ld} |f(\og)|^p \tr_{_{1/n}} (\og) \leq
C \tau \int_{B(e,\al)} |f(x)|^p\, d\sa(x),$$ where    $\tau$ is
defined by
$$ \tau=\max_{x\in\sph} \#\(\Ld \cap B_{\rho} ( x, \f 1n)\).$$
\end{cor}

\begin{cor}
If    $\be \ge 1$, $1\leq p<\infty$, and $\Ld$ is a maximal $(\f
\da n,\rho)$-separable subset of $B(e,\al)$ with $\da\in (0,\f
1{4C_1}]$ and  $C_1$  the same as  in Equation ($\ref{1-6}$), then
for all $f\in\Pi_n^d$, we have
\begin{align*} \int_{B(e,\al)} |f(x)|^p \, d\sa(x)& \sim \(\sum_{\og \in \Ld} \(\max_{x\in
B_\rho ( \og, \be \da/ n)} |f(x)|^p\) \tr_{_{\da/n}}
(\og)\)^{\f1p}\\
&\sim \(\sum_{\og \in \Ld} \(\min_{x\in B_\rho ( \og,\,  \be\da
/n)} |f(x)|^p\)\tr_{_{\da/n}} (\og)\)^{\f1p},
\end{align*} where the constants of equivalence are independent of
$f$, $n$, $\al$ and $\{\og\}_{\og\in\Ld}$.  \end{cor}

{\it Remark 1.5.}\  \    In one dimensional case,   the following
large sieve inequality was  proved by  Golinskii, Lubinsky, and
Nevai [GLN]:
\begin{equation} \label{1-8-c} \sum_{k=1}^m |P(\al_j)|^p \va(\al_j) \leq C \tau
\int_a^b  |P(\ta)|^p\, d\ta\end{equation} with $C$ independent of
$m$, $n$, $P$, $p$, $a$, $b$ ,$\{\al_j\}$.  Here
   $P$ is a trigonometric polynomial of degree $\leq
n$, \footnote{In [GLN],  $P$ could be a `` generalized
trigonometric polynomial'', not just an ordinary trigonometric
polynomial. }
$$\va(\ta)= \f 1{ pn+1} \( \Bl| \sin \(\f {\ta-a}2\)\sin \( \f
{\ta-b} 2\)\Br|+\(\f {b-a}{pn+1}\)^2\)^{1/2},$$
 while
$$ 0\leq a \leq \al_1<\al_2<\cdots <\al_m \leq b \leq
2\pi,$$

$$\tau =\max_{\ta\in[a, b]} \#\(\{j:\   \  \al_j
\in[\ta-\va(\ta),\ta+\va(\ta)]\}\), $$ $0<p<\infty$  and $m\ge 1$.
 A
version of ($\ref{1-8-c}$), which has the correct form for all
choices of $[a, b]$ -whether $b-a$ is very small or close to
$2\pi$, was established in a recent paper [KL2] by Kobindarajah
and Lubinsky. Note that if  we identify the interval $[a,b]$ with
the arc $B(e, \al)$ centered at $ e=(\cos\f {a+b}2, \sin \f
{a+b}2)$ and of radius $\al= \f {b-a}2$, we would have
  $$ \va(\ta)\sim \tr_{_{1/n}} (z)\sim \Bl|B_\rho(z, \f 1n)\Br|,$$
  provided that  $b-a$ is not too close to $2\pi$,
   where $\ta\in [a,b]$ and $z=(\cos\ta,\sin\ta)$.
This means that  Corollary 1.3 can be considered  as a
higher-dimensional analogue of the  large sieve inequality
($\ref{1-8-c}$). While we believe Corollary 1.3 remains true for
$0<p<1$ as well, we are unable to prove it.

{\it Remark 1.6.}  \    It was shown by Mhaskar [M1] that given an
arbitrary set $\Ld$ of points in $B(e,\al)$ satisfying the mesh
norm condition $\dmax_{\xi\in\sph} \dmin _{\og\in\Ld}d( \og,\xi)
\leq c\al$, there exists nonnegative weights $\ld_\og$,
$\og\in\Ld$ such that ($\ref{1-1}$) holds for every $f\in\Pi_n^d$.
Here  the constant $c$   is independent of $f$ and $\al$, but
depends on the degree $n$. Here  we wish to compare our result
with that of Mhaskar [M1]. First, our results are uniform in the
degree $n$ and the radius $\al$ (when $\al$ is not very close to
$\pi$), while his result is not uniform in $n$. Indeed, in our
opinion, uniformity in  the degree is of crucial importance, as
can be seen in many known  work (see, for instance, [GLN, D2,
NPW2]).
 Second, our result ( Corollary 1.2) shows  that the minimum  number of
nodes  required in  a positive cubature formula of degree $n$ on
$B(e,\al)$ is comparable to the dimension of the space $\Pi_n^d$
as $n\to\infty$, while his result does not.  Third, we have a
sharp estimate $\ld_\og \sim \tr_{_{\da/n}} (\og)$ on  the weights
of the cubature formula (\ref{1-1}), while only nonnegativity of
the cubature  weights were shown in [M1]. Fourth, as mentioned in
Remark 1.5, our result can be considered as a higher dimensional
analogue of the large sieve inequality of Golinskii, Lubinsky, and
Nevai. Fifth, in our work we made special efforts to treat the
boundary and the center of the spherical cap, while Mhaskar [M1]
didn't. In our opinion,  good cubature formulas and MZ
inequalities on a spherical cap $B(e,\al)$  couldn't be obtained
without taking into consideration of the boundary of $B(e,\al)$.
Indeed, our opinion is supported by many known results on a finite
interval $[a, b]$. Last but not least, as is demonstrated in
Section 5 (Theorem 5.1), our method can be used  to  obtain  MZ
inequalities with doubling weights on spherical caps.

The paper is organized as follows.   In Section 2,  we show  two
technical lemmas concerning  the properties of the metric $\rho$
in the case when $\al \in (0,\f12]$. After that, we prove the main
results, Theorem 1.1 and Corollaries 1.2--1.4, for the case $\al
\in (0,\f12]$  in Section 3. The proofs of the main results for
the remaining case $\al \in (\f 12,\pi)$ can be deduced from the
case $\al\in (0,\f 12]$. This is done in Section 4. Finally, in
Section 5, we discuss briefly how to establish similar results for
spherical collars and for spherical caps with   doubling weights.

\section{Two basic lemmas}

In this section, we establish some basic  facts concerning the
metric $\rho\equiv \rho_{_{B(e,\al)}}$ defined by ($\ref{1-3}$) in
the case when $\al\in (0,\f12]$. These facts will be needed in
later sections. We shall begin with the simple case $d=1$, where
$\mathbb{S}^1$ is the unit circle, identified as $\mathbb{R}/ 2\pi
\mathbb{Z}$.

 Let $\al\in (0,\f 12)$.  For $x_1, x_2\in [-\al,\al]$, we define
\begin{equation}\label{2-1-a}\rho_1(x_1, x_2)\equiv \rho_{_{[-\al,\al]}}(x_1,x_2):=
\frac 1\al \sqrt{|x_1-x_2|^2+\al
\,\Big|b_{x_1}^{1/2}-b_{x_2}^{1/2} \Big|^2}, \end{equation}
 where $b_x\equiv b_{_{x,\  [-\al,\al]}}$ denotes the shortest distance
 from $x\in [-\al,\al]$ to the boundary of $[-\al,\al]$, that is,
$$b_{x}\equiv b_{x, [-\al,\al]}:=\min\{|x+\al|,|x-\al|\}.$$   Clearly, $\rho_1$ is the one-dimensional
analog  of the metric $\rho\equiv \rho_{_{B(e,\al)}}$ defined by
($\ref{1-3}$). It turns out ( see Lemma 2.1 below) that $\rho_1$
is equivalent to  two other metrics $\rho_2$ and $\rho_3$ on
$[-\al,\al]$, whose definitions are given as follows:
 for $ x_1 =\arcsin ((\sin\al) \cos t_1)$ and $ x_2 =\arcsin
((\sin\al) \cos t_2)$ with $t_1, t_2\in [0,\pi]$,
\begin{align}
\rho_2 ( x_1, x_2)&= \f 1\al \sqrt{ |x_1-x_2|^2 + \Bl| \sqrt{
\al^2-x_1^2} - \sqrt{ \al^2-x_2^2}\Br|^2},\label{2-3}\\
\rho_3( x_1, x_2)&=|t_1-t_2|.\label{2-4}\end{align} For $x\in
[-\al,\al]$ and $r\in (0, 1)$, we write $B_{\rho_i} (x, r) :=\{
y\in[-\al,\al]:\   \ \rho_i (x, y) \leq r\},$  $i=1,2,3$.

Now our first lemma can be stated as follows:

\begin{lem}
Let $\al\in (0,\f12]$ and  let $\rho_1$, $\rho_2$ and $\rho_3$ be
defined by ($\ref{2-1-a}$), ($\ref{2-3}$) and ($\ref{2-4}$),
respectively. Then the following statements hold true:

{\rm (i)}\   \   For any $x_1, x_2\in [-\al,\al]$,
\begin{equation}\label{2-5} \rho_1(x_1, x_2)\sim \rho_2(x_1, x_2)\sim \rho_3(x_1,
x_2),\end{equation} where the constants of equivalence are
independent of $\al$, $x_1$ and $x_2$.

{\rm (ii)} \  \   For any  $x\in [-\al,\al]$ and $r\in (0, 1)$, we
have  $$B_{\rho_1} (x, r) \subset [ x- \al r, x+\al r]\cap [-\al,
\al]$$ and
\begin{equation}\label{2-6}
|B_{\rho_i} (x, r)|\sim \al  \( r^2 + r \sqrt{1-(x/\al)^2}\),\   \
i=1,2,3,
\end{equation} where the constant of equivalence is
independent of $x$, $r$ and $\al$.

{\rm (iii)}\   \  For any $x_1, x_2\in [-\al,\al]$ and $r>0$,
\begin{equation}\label{2-6-a}|B_{\rho_1} (x_1, r)| \leq C \( 1+ \f{\rho_1(x_1, x_2)} r\)
|B_{\rho_1} (x_2, r)|,\end{equation} where $C>0$ is independent of
$x_1, x_2, r$ and $\al$.

{\rm (iv)}\   \   If $r\in (0, 1)$, $\be\ge 1$ and  $A$ is an $(r,
\rho_1)$-separable subset of $[-\al,\al]$, then
\begin{equation}\label{2-7-a} \sup_{x\in[-\al,\al]} \sum _{ \xi\in A} \chi_{_{
B_{\rho_1} (\xi, \be r)}} (x) \leq C\be^3,\end{equation} where $C$
is an absolute constant.

\end{lem}

\begin{proof}\   \ (i)\   \
First, we show the equivalence
\begin{equation}\label{2-7}\rho_1(x_1, x_2)\sim \rho_2(x_1,
x_2).\end{equation} To show this, we  start with the case
$x_1x_2\ge 0$. Without loss of generality, we may assume in this
case that $x_1, x_2\in [0,\al]$  (otherwise, consider $-x_1$ and $
-x_2$).  Then $ b_{x_i}=\al-x_i,\    \  i=1,2,$ and thus, by the
definition, we obtain
\begin{align}
\rho_1(x_1,x_2)& \sim\f {|x_1-x_2|} \al+\f {|x_1-x_2|}{\sqrt{\al}
( \sqrt{\al-x_1}
+\sqrt{\al-x_2})},\label{2-8}\\[4mm]
\rho_2(x_1,x_2)& \sim \f {|x_1-x_2|} \al+\f {|x_1-x_2|}{\sqrt{\al}
( \sqrt{\al-x_1} +\sqrt{\al-x_2})}\f {x_1+x_2}
{\al}.\label{2-9}\end{align} If $x_1+x_2\ge \f \al4$ then
$x_1+x_2\sim \al$, so  comparison  ($\ref{2-8}$) with
($\ref{2-9}$) gives $\rho_1(x_1, x_2)\sim \rho_2(x_1, x_2).$
However, on the other hand, if $x_1+x_2\leq \f \al 4$ then $
\sqrt{\al-x_1}+\sqrt{\al-x_2} \sim\sqrt{\al},$ so by ($\ref{2-8}$)
and ($\ref{2-9}$), we deduce $\rho_1(x_1, x_2)\sim \rho_2(x_1,
x_2)\sim \f{|x_1-x_2|}\al.$ This proves ($\ref{2-7}$) in the case
$x_1x_2 \ge 0$.

Equation ($\ref{2-7}$) for the case $x_1x_2<0$ follows from the
case $x_1x_2\ge 0$.  In fact, by the already proven case
$x_1x_2\ge 0$, we deduce that if $x_1x_2 <0$ then  $\rho_1(x_1,
-x_2)\sim \rho_2(x_1, -x_2)$, which together with  ($\ref{2-1-a}$)
and ($\ref{2-3}$) implies \begin{equation}\label{2-10}\f
{|x_1+x_2|} \al + \f { |\sqrt{\al^2-x_1^2} -\sqrt{\al^2-x_2^2}|}
\al \sim \f {|x_1+x_2|} \al + \f
{|b_{x_1}^{\f12}-b_{x_2}^{\f12}|}{\sqrt{\al}}.\end{equation} Note,
on the other hand,  that if $x_1\cdot x_2\leq 0$ then $
|x_1+x_2|\leq |x_1-x_2|.$  Thus,  using ($\ref{2-10}$), we
conclude that for $x_1, x_2\in[-\al,\al]$ with $x_1x_2<0$,
$$\f {|x_1-x_2|} \al + \f { |\sqrt{\al^2-x_1^2}
-\sqrt{\al^2-x_2^2}|} \al \sim \f {|x_1-x_2|} \al + \f
{|b_{x_1}^{\f12}-b_{x_2}^{\f12}|}{\sqrt{\al}},$$ which implies
$\rho_1(x_1, x_2)\sim \rho_2(x_1, x_2)$, and therefore  completes
the proof of ($\ref{2-7}$).

Next, we show
\begin{equation}\label{2-11}
\rho_2(x_1, x_2) \sim \rho_3(x_1, x_2).\end{equation} To this end,
 we set, for $t\in [0,\pi]$,
\begin{equation*}\label{2-1-1} g(t) = \arcsin ((\sin\al)\cos
t)   \    \  \text{and}\   \  h(t) = \sqrt{ \al^2 -
(g(t))^2}.\end{equation*} Since $\al \in (0,\f 12]$, it is easy to
verify that for $t\in [0,\pi]$
\begin{equation}\label{2-13-a}
g(t)\sim h'(t)\sim \al \cos t \    \    \text{and}\    \    \  \
h(t)\sim -g'(t) \sim \al \sin t.\end{equation}
 Now we assume that $x_1=g(t_1)$ and
$x_2=g(t_2)$ with $t_1, \  t_2 \in [0,\pi]$. Then by
($\ref{2-3}$),
\begin{align}
\rho_2(x_1, x_2) &\sim  \f 1\al \Bl[ |g( t_1)-g(t_2)| +
|h(t_1)-h(t_2)|
\Br]\notag\\
& = \f1\al \Bl| \int_I g'(t)\, dt \Br| +\f1\al \Bl| \int_I h'(t)\,
dt \Br|,\label{2-14}\end{align} where $I=[t_1, t_2]$ or $[t_2,
t_1]$. By ($\ref{2-13-a}$), we obtain
$$ \rho_2(x_1, x_2) \leq C |I|=C |t_1-t_2|=C \rho_3 (x_1, x_2).$$
To show the converse inequality
$$ \rho_2(x_1, x_2) \ge C |t_1-t_2|,$$
 we note that if $|t_1-t_2|\leq \f \pi6$ then by the mean value
theorem for integrals, we obtain, for some $\xi_1, \ \xi_2 \in I$,
\begin{align*}
\rho_2(x_1, x_2) &=\f {|I|}\al \( |g'(\xi_1)|+|h'(\xi_2)|\)\\
&\ge C |I| (|\sin \xi_1| +|\cos \xi_2| )\    \      \   \    \ \
\text{(by ($\ref{2-13-a}$))}\\
&\ge C |I|\(\sin^2\xi_1+\cos^2\xi_1 - | \cos \xi_2-\cos
\xi_1|\)\\
& \ge C |I|\( 1-|\xi_1-\xi_2|\)\ge C(1-\f \pi 6)|I|.\end{align*}
On the other hand,  if $ |t_1-t_2|\ge \f \pi6$,  then by
($\ref{2-13-a}$) it follows that
$$\rho_2(x_1, x_2)\ge \f 1\al |\int_I g'(t)\, dt | \ge  C \int_I \sin
t\, dt\ge C |I|.$$  This completes the proof of ($\ref{2-11}$).

(ii)\   \ Since by the definition, for all $x_1,
x_2\in[-\al,\al]$,
$$|x_1-x_2|\leq  \al \rho_1(x_1, x_2),$$
it follows that  $B_{\rho_1} (x, r)\subset [x-\al r, x+\al
r]\bigcap [-\al,\al]$.  Thus,  by ($\ref{2-5}$), it remains  to
show
\begin{equation}\label{2-15-a}
|B_{\rho_3} (x, r)|\sim \al  \( r^2 + r \sqrt{1-(x/\al)^2}\).
\end{equation}
Again, we set $g(t)=\arcsin ( (\sin\al)\cos t)$. Given
$x\in[-\al,\al]$, we shall use the natation $t_x$ to denote the
unique solution in $[0,\pi]$ to the equation $g(t)=x$. Then we
have
\begin{equation} \label{2-15}B_{\rho_3} (x, r) =\Bl\{ g(t):\ \ t\in
[0,\pi]\    \ \text{and}\
  \   |t-t_x|\leq r\Br\}.\end{equation}

For the proof of ($\ref{2-15-a}$), we start with the case $ x\in [
0,\al]$. In this case,  $t_x\in [0,\f \pi2]$, and therefore
setting  $ \ga =\max \{ 0, t_x-r\}$, we obtain
\begin{align}
| B_{\rho_3} ( x, r) |& = g(\ga) - g( t_x+r)  \sim \al
\int_\ga^{t_x+r} \sin u \, du\notag\\
&\sim \al \int_{t_x+\f r2} ^{t_x+r} \sin u \, du \sim \al r
(t_x+r),\label{2-16}\end{align} where in the first ``$\sim$'' we
have used ($\ref{2-13-a}$), while in the second ``$\sim$'' we have
used  the doubling property of the weight function  $|\sin t|$ (
see, for instance, [MT2]). On the other hand, by ($\ref{2-3}$) and
($\ref{2-5}$), we have
$$ t_x= \rho_3 (x, \al)\sim \rho_2 (x, \al) \sim \sqrt { 1- \f
x\al}\sim \sqrt { 1- \(\f x\al\)^2},$$  which combined with
($\ref{2-16}$) yields the desired equation ($\ref{2-15-a}$) in the
case $x\in [0,\al]$.

We conclude the proof of ($\ref{2-15-a}$) by showing that  the
case $x\in [-\al,0]$ follows from the already proven case $x\in
[0,\al]$. In fact, since $g(\pi-t)=-g(t)$, we have, for $x\in
[-\al, 0]$,
 \begin{align*}
 B_{\rho_3}
(x, r)& =\{ g(t):\   \  t\in [0,\pi]\   \  \text{and}\
  \   |t-t_x|\leq r\}\\
&=\{- g(u):\   \  u\in [0,\pi]\   \  \text{and}\
  \   |u-t_{-x}|\leq r\}= -B_{\rho_3} (-x, r).\end{align*}
  Thus,  by the already proven case $x\in [0,\al]$, we deduce that  for
  $x\in [-\al,0]$
  $$ | B_{\rho_3}
(x, r)|=| B_{\rho_3} (-x, r)| \sim \al \(  r^2 + r\sqrt{
1-(x/\al)^2}\),$$ which gives  ($\ref{2-15-a}$) in this case.

(iii)\   \ Note that   for $x_1, x_2\in[-\al,\al]$,
$$\Bl| \f { \al ( r^2+r\sqrt{ 1- (x_1/\al)^2})}{ \al ( r^2+r\sqrt{ 1-
(x_2/\al)^2})}-1\Br|
 = \f{|\sqrt{\al^2-x_1^2}-\sqrt{\al^2-x_2^2}|}{ \al r+\sqrt{
 \al^2-x_2^2}}
 \leq \f{\rho_2 (x_1, x_2)}r.$$
The desired inequality  ($\ref{2-6-a}$) then follows by
($\ref{2-6}$) and ($\ref{2-5}$).

(iv)\    \ Let $\be\ge 1$ and let  $A$ be  an $(r,
\rho_1)$-separable subset of $[-\al,\al]$. Then by the definition
of $(r,\rho_1)$-separable, it follows that for any  $x\in
[-\al,\al]$, \begin{equation}\label{2-18-a} \sum_{\xi\in A\cap
B_{\rho_1} (x,\be r)} | B_{\rho_1} (\xi,\f r4)|\leq \bl|
B_{\rho_1} (x,(\be+\f 14)r)\br|.\end{equation} However, on the
other hand, by ($\ref{2-6-a}$), we note that for any $\xi\in
B_{\rho_1} (x, \be r)$, \begin{equation}\label{2-19-a} (C\be)^{-1}
|B_{\rho_1} (x, \f r4)| \leq |B_{\rho_1} (\xi, \f r4)|\leq C\be
|B_{\rho_1} (x, \f r4)|.\end{equation} Thus, combining
($\ref{2-18-a}$) with ($\ref{2-19-a}$),  we obtain
$$\sum _{\xi\in A} \chi_{_{B_{\rho_1} (\xi, \be r)}} (x) =\#
\( A \bigcap B_{\rho_1} (x, \be r)\) \leq C \be
 \f {\bl| B_{\rho_1} (x,(\be+\f14)r)\br|}{\bl| B_{\rho_1} (x,\f
r4)\br|}\leq C \be ^3,$$ which proves ($\ref{2-7-a}$).

The  proof of Lemma 2.1 is complete.
\end{proof}

Now we turn to the case $d\ge 2$. Recall that
$\rho=\rho_{_{B(e,\al)}}$ is the metric on the spherical cap
$B(e,\al)$ defined by ($\ref{1-3}$), and
$\rho_1=\rho_{[-\al,\al]}$ is the metric on $[-\al,\al]$ defined
by ($\ref{2-1-a}$). We need to introduce two more metrics $\rho_4$
and $\rho_5$ on $B(e,\al)$. To this end, we set, for $e\in\sph$,
$$ \mathbb{S}_e^{d-1}: =\{ y\in\sph:\   \  y\cdot e =0\}.$$
For $x=e\cos \ta +\xi \sin\ta$ and $y= e\cos t+\eta\sin t$ with
$\xi,\eta\in \mathbb{S}_e^{d-1}$ and $\ta, t\in [0,\al]$, we
define
\begin{equation}\label{2-20-a}
\rho_4(x,y):=\max\Bl\{ \rho_1 (\ta, t),\
d(\xi,\eta)\Br\},\end{equation} \begin{equation}\label{2-21-a}
\rho_5(x, y): = \f 1{\sin\al} \sqrt{ |\xi \sin\ta-\eta\sin t|^2 +
\bl| \sqrt{\sin^2\al - \sin^2 \ta}-\sqrt{\sin^2\al - \sin^2
t}\br|^2}.\end{equation} Recall that for   $x\in B(x,r)$,  $r>0$
and a metric $\widetilde{\rho}$ on $B(e,\al)$,
$$ B_{\widetilde{\rho}} (x, r)=\Bl\{y\in B(e,\al):\   \   \wt{\rho}
(x, y) \leq r\Br\}.$$

\begin{lem}
Let   $\va\in (0,1)$, $e\in\sph$ and $\al\in (0,\f12]$. Then the
following statements hold true:

 (i)  \begin{align} &\rho(x,y)\sim \rho_5(x,y), \ \ \text{for
all $x, y \in B(e,\al)$;}\label{2-22-a}\\
&\rho(x,y) \sim \rho_4(x,y), \   \   \   \text{if $x, y\in
B(e; \va\al, \al)$ };  \label{2-23-a} \\
 &\rho(x,y) \sim \f 1\al d(x,y),\    \  \text{if $x, y\in B(e,
(1-\va) \al)$},\label{2-24-a}\end{align}where the constants of
equivalence are independent of $x, y$ and $\al$, but may depend on
$\va$ when $\va$ is small.

(ii) If $x\in B(e,\al)$, then for any
$r>0$,\begin{equation}\label{2-25-a} B_{\rho} (x, r) \subset
B(x,\al r),\ \   B_{\rho_4} (x, r) \subset B(x,3\al
r);\end{equation} if $x \in B(e;\va\al, \al)$, then for any $r>0$,
\begin{equation}\label{2-26-a} B_{\rho_4} (x, C_2^{-1} r)\subset  B_{\rho}
(x, r)\subset  B_{\rho_4} (x, C_2r);\end{equation} if  $x\in B(e,
(1-\va) \al)$ then for any $r>0$, \begin{equation}\label{2-27-a} B
(x, C_2^{-1}\al r)\subset B_{\rho} (x, r)\subset  B (x, \al
r),\end{equation}where $C_2$ is independent of $r$, $x$ and $\al$,
but depends on $\va$ when $\va$ is small.

(iii)\   \  For  any  $x\in B(e,\al)$ and $r\in (0,1)$,
\begin{equation}\label{2-28-a} | B_\rho(x,r) | \sim \al^{d} \( r^{d+1} + r^{d}
\sqrt{\f{ b_{x}}\al}\),\end{equation} where $b_x\equiv b_{x,
B(e,\al)}$ is defined by ($\ref{1-4}$), and the constant of
equivalence depends only on $d$.

(iv) \    \  For any $x, y\in B(e,\al)$ and $r>0$,
\begin{equation}\label{2-29-a} | B_\rho(x,r) | \leq C \( 1+ \f {\rho(x,y)} r\) |
B_\rho(y,r)|,\end{equation}where $C>0$ depends only on $d$.

(v) \    \  Suppose that  $r\in (0,1)$, $\be \ge 1$ and  $\Ld$ is
an $(r, \rho)$-separable subset of $B(e,\al)$. Then we have
\begin{equation}\label{2-30-a} \max_{x\in B(e,\al)} \sum _{\og\in \Ld}
\chi_{_{B_{\rho} (\og, \be r)}}(x) \leq
C\be^{d+2},\end{equation}where $C>0$ depends only on $d$.

\end{lem}

\begin{proof} \   \  (i) \   \   Let
$x=\xi\sin\ta+e\cos\ta$ and $y=\eta\sin t+e\cos t$ with $\xi,
\eta\in\mathbb{S}_e^{d-1}$ and $\ta,t \in [0,\al]$. We start with
the proof of ($\ref{2-22-a}$). We first note that
 \begin{equation}\label{2-31-a} 4\sin^2\(\f{d(x,y)} 2\) =
 |x-y|^2 =4\sin^2\f{\ta-t}2 +(\sin\ta\sin t)
|\xi-\eta|^2,\end{equation} which implies
\begin{equation}\label{2-32-a}
d(x, y) \sim |\ta-t| + |\eta-\xi|\sqrt{\ta\cdot  t}.\end{equation}
Since $\al\in (0,\f12]$, it follows by a straightforward
calculation  that \begin{equation}\label{2-33-a}|\xi\sin\ta -\eta
\sin t|\sim |\ta-t| + |\eta-\xi|\sqrt{\ta\cdot t} \sim
d(x,y),\end{equation} and
\begin{equation}\label{2-34-a}
\bl|\sqrt{\sin^2\al -\sin^2\ta} -\sqrt{\sin^2\al -\sin^2 t}\br|
\sim\f {(\ta+t) |\sqrt{b_x} -\sqrt{b_y}|}{\sqrt{\al}},
\end{equation} where $b_x=\al-\ta$ and $b_y=\al-t$.
On the other hand, note, however, that $\ta+t \in [0, 2\al]$ and
that
$$ \sqrt{\al} |\sqrt{b_x}-\sqrt{b_y}| =\f { \sqrt{\al}
|\ta-t|}{\sqrt{\al-\ta} +\sqrt{\al-t}}\sim |\ta-t|\leq C d(x,y),$$
provided  $\ta+t< \f \al2$. This means that
\begin{equation}\label{2-35-a}d(x, y) + \f {(\ta+t) |\sqrt{b_x}
-\sqrt{b_y}|}{\sqrt{\al}}\sim d(x, y) + \sqrt{\al}
|\sqrt{b_x}-\sqrt{b_y}|.\end{equation}
 Therefore, combining ($\ref{2-33-a}$)--($\ref{2-35-a}$), we
deduce the desired equivalence ($\ref{2-22-a}$).

Next, we show ($\ref{2-23-a}$) in the case when $\ta, t \in
[\va\al, \al]$. In fact, we have
\begin{align*}
\rho_4(x,y) &\sim |\xi-\eta| +\Bl[\f
1\al |\ta-t| +\f 1{\sqrt{\al}} |\sqrt{\al-t}-\sqrt{\al-\ta}|\Br]\\
&\sim \f 1\al d(x,y) +\f 1{\sqrt{\al}}
|\sqrt{\al-t}-\sqrt{\al-\ta}| \sim \rho(x,y),\notag\end{align*}
where in the first ``$\sim$'' we have used ($\ref{2-1-a}$) and
($\ref{2-20-a}$), in the second ``$\sim$'' we have used
($\ref{2-32-a}$) and the fact that $\ta, t\in[\va \al,\al]$, and
the final ``$\sim$'' follows by ($\ref{1-3}$). This proves the
desired equation ($\ref{2-23-a}$).

 Finally, we note that   ($\ref{2-24-a}$) for $\ta, t \in [0, (1-\va) \al]$ is a simple
 consequence of the definition ($\ref{1-3}$) and the following equation:
\begin{align*}\f 1{\sqrt{\al}} |\sqrt{b_x} -\sqrt{b_y}|
=\f 1{\sqrt{\al}} \f{|\ta-t| }{ \sqrt{\al-t} +\sqrt{\al-\ta}} \sim
\f1 \al |\ta-t|\leq \al^{-1} d(x,y).\end{align*}

(ii)\   \  It follows   by ($\ref{1-3}$), ($\ref{2-20-a}$) and
($\ref{2-31-a}$) that
\begin{equation}\label{2-36-a} d(x,y)\leq \min\{ \al \rho(x,y), 3\al
\rho_4(x,y)\},\    \   \  x, y\in B(e,\al),  \end{equation} which
implies ($\ref{2-25-a}$).  Thus, it remains to show
($\ref{2-26-a}$) and ($\ref{2-27-a}$).  By the definition, it's
easily seen that for all $u, v\in [0,\al]$ and $y, z\in B(e,\al)$,
$$\rho_1(u, v) \leq 3,\    \    \    \   \max\{ \rho(y, z),\  \rho_4(y, z)\} \leq \pi.$$ Thus, without loss of generality, we may
assume that $r\in (0, \f \va {6}]$.   Then taking into account
($\ref{2-36-a}$), we deduce that  for $x\in B(e; \va \al, \al)$,
$$ B_{\rho} (x, r)\bigcup B_{\rho_4} (x, r) \subset B (e; \f
{\va\al} 2, \al),  $$ which together with ($\ref{2-23-a}$) implies
($\ref{2-26-a}$). Finally, Equation ($\ref{2-27-a}$) follows by
($\ref{2-24-a}$) and ($\ref{2-25-a}$).

(iii)\   \  We start with the case $\al/6\leq \ta: =d(x,e)\leq
\al$. In this case, by ($\ref{2-26-a}$), it is sufficient to show
that for $r\in (0,\f 1{12})$,
\begin{equation}\label{2-37-a} |B_{\rho_4} (x,r) | \sim \al^{d} \( r^{d+1} +r^{d}
\sqrt{b_x /\al}\).\end{equation} Notice that by Lemma 2.1 (ii),
$B_{\rho_1}(\ta,r)\subset [\f 1{12}\al,\al]$.  Thus $$ |B_{\rho_4}
(x,r)|\sim r^{d-1}\int_{B_{\rho_1}(\ta,r)} \sin^{d-1} t\, dt \sim
(\al r)^{d-1} |B_{\rho_1} ( \ta, r)|.$$ This last equation
together with ($\ref{2-6}$) implies  ($\ref{2-37-a}$) and hence
($\ref{2-28-a}$) in the case when $\ta =d(x, e) \ge \f 16 \al$.

Finally, we note that   ($\ref{2-28-a}$) for  the case  $0\leq
\ta=d(x,e)\leq \f \al 6$ follows directly from ($\ref{2-27-a}$).

(iv)\   \   Inequality ($\ref{2-29-a}$) is a simple consequence of
($\ref{2-28-a}$) and the following equation:
$$ \Bl| \f { \sqrt{\al} r+ \sqrt{b_x}} { \sqrt{\al} r+
\sqrt{b_y}}-1\Br|=\f {|\sqrt{b_x}-\sqrt{b_y}|} { \sqrt{\al} r
+\sqrt{b_y}}\leq \f {\sqrt{\al} \rho(x,y)} { \sqrt{\al}
r+\sqrt{b_y}}\leq \f { \rho(x,y)} r.$$

(v) \   \ ($\ref{2-30-a}$) follows by ($\ref{2-28-a}$),
($\ref{2-29-a}$) and the standard volume comparison method. Since
the proof is almost identical to that of Lemma 2.1 (iv), we omit
the details.

 This completes the proof of
Lemma 2.2.

\end{proof}

\section{Proofs of the main results for $\al \in (0,\f12]$}

The proofs of Theorem 1.1 and Corollaries 1.2--1.4 in the case
when $\al\in (0,\f12]$ are  based on a series of lemmas. To state
these lemmas, we need to introduce several notations.
  We say a weight function $W$ on $\sph$ is a doubling weight if
there exists a constant $L$, called doubling constant, such that
for all $x\in\sph$ and $r\in (0,\pi)$,
$$ W\(B(x,2r)\) \leq L  W\(B(x,r)\),$$
   where and elsewhere, we write, for
a subset  $E$ of $\sph$,
$$ W(E) =\int_E W(y)\, d\sa(y).$$
As usual, we identify the unit circle $\mathbb{S}^1$ with
$\mathbb{R}/2\pi \mathbb{Z}$. Thus, $\Pi_n^1\equiv
\Pi_n(\mathbb{S}^1)$ denotes the space of all trigonometric
polynomials of degree at most $n$ on $\mathbb{R}$. Associated with
a  function $f$ on $[-\al,\al]$, we define
$$ f_\al (t)=f\bl(\arcsin ((\sin \al) \cos t)\br),\   \   t\in [-\pi,\pi],$$
and associated with a weight function $W$ on $[-\al,\al]$, we
define
$$ W_{n,\al} ( t) := n \int_{t-\f 1n}^{t+\f 1n} W_\al (\ta)\,
d\ta,\    \  n=1,2,\cdots.$$

 Our first
lemma is due to T. Erd\'elyi [Er2, Theorems 1.3 and 2.1]:

\begin{lem} \mbox{\rm ([Er2]).}\    \
Let $p\in [1,\infty)$ and $\al\in (0,\f12]$. Suppose $W$ is a
weight function on $[-\al,\al]$ such that $W_\al $ is a doubling
weight on $\mathbb{S}^1$. Then for all $T\in\Pi_n^1$, we have
\begin{equation}\label{3-1} \int_{-\al}^\al |T'(t) | ^p W(t) \( \f \al n + \sqrt{
\al^2 -t^2}\)^p \, dt\leq C n^p \int_{-\al} ^\al |T(t)|^p W(t)\,
dt,\end{equation} and
\begin{equation}\label{3-2-a}
  \int_{-\pi}^\pi |T_\al (t) |^p W_\al (t) |\sin t
|\, dt\\
 \sim \int_{-\pi}^\pi |T_\al (t) |^p W_{n,\al} (t) |\sin t |\,
dt,\end{equation}
  where  the constant $C$  and the constant of equivalence   depend only on $p$
and the doubling constant of $W_\al$.
\end{lem}

It was pointed out in [Er2] that $W_\al$   is a doubling weight if
and only if  $W(\al\cos t)$ is a doubling weight.
 In the unweighted case, ($\ref{3-1}$) for all $0<p<\infty$ was
proved by Lubinsky [L] ( see also the paper [KL1] by Kobindarajah
and Lubinsky ). For relevant results concerning doubling weights,
one may consult \mbox{[D1, Er1,Er2, MT1, MT2, MT3]}.

To state our next lemma, we recall that
$\rho_1=\rho_{_{[-\al,\al]}}$ is the metric on $[-\al,\al]$
defined by ($\ref{2-1-a}$).

\begin{lem} Let $\al \in (0, \f 12]$, $\be \ge 1$, $1\leq p<\infty$ and
 $\da\in (0,1)$. Let  $W$ be a weight function on
$[-\al,\al]$ such that $W_\al$ is a doubling weight on
$\mathbb{S}^1$.
  Suppose that $n$ is  a positive integer and   $\{\xi_j\}_{j=1}^{m_n}$  is a
 $(\f\da n, \rho_1)$-separable  subset of $
[-\al,\al]$.    Then
 for any $T\in\Pi_n^1$,
we have \begin{equation*} \sum _{i=1}^{m_n} \(\max_{x,\, y\in
B_{\rho_1} (\xi_i, \f {\be \da} n)} |T(x)-T(y)|^p\)
\int_{B_{\rho_1} (\xi_i, \f \da n)} W(t)\, dt \leq (C\be^2\da)^p
\int_{-\al} ^\al |T(x)|^p W(x)\, dx,\end{equation*} where $C>0$
depends only on $p$ and  the doubling constant of $W_\al$.
\end{lem}
\begin{proof}
 As in the proof of Lemma 2.1, we set
$$ g(t)\equiv g(t,\al)= \arcsin ((\sin\al) \cos t),\    \  t\in [-\pi,\pi].$$
The proof is based on Lemma 3.1.  Let  $n_1=[n/\da]$,
$T\in\Pi^1_n$, and  $1\leq p<\infty$.
  Suppose
  that
$$\{\xi_j =g(t_j):\    \  t_j\in[0,\pi],\   \  j=1,2,\cdots m_n\}$$
is $(\f \da n, \rho_1)$-separable in $[-\al,\al]$.  Then  by Lemma
2.1 (i),   there exists an absolute constant $\ga \ge 1$ such that
\begin{equation}\label{3-4-b}
\min_{1\leq i\neq j\leq m_n} |t_i-t_j| \ge \f {\da} {\ga
n},\end{equation} and for all $1\leq j \leq m_n$ and  $ r>0$,
\begin{equation}\label{3-5-b} B_{\rho_1} (\xi_j, r) \subset \Bl\{ g(t):\    \  t\in
[t_j-\ga r, t_j+\ga r]\cap [0,\pi]\Br\}.\end{equation} It follows
that for a fixed $j\in [1, m_n]$ and any $\be \ge 1$,
\begin{align*}
\( &\max_{x, y \in B_{\rho_1} ( \xi_j, \f{\be\da} n)}
|T(x)-T(y)|^p\)
\(\int_{B_{\rho_1}(\xi_j, \f \da n)}  W(\xi)\, d\xi \)\\
&\leq C \al^{p+1}  \( \int_{t_j-\f {\ga \be \da}n }^{t_{j}+\f {\ga
\be \da}n} |T'(g(t))| ^p |\sin t| \, dt \)\( \int_{t_j-\f {\ga
\da}n}^{t_{j}+\f {\ga \da}n} W_\al (t) |\sin t|\, dt
\) \(\int_{t_{j}-\f {\ga \be \da}n}^{t_{j}+\f {\ga \be \da}n}|\sin t| \, dt \)^{p-1}\\
&\leq C   \(\f{\ga \be \da}{n}\)^{p-1} \al^{p+1} \min_{t\in
[t_{j}-\f {\ga \be \da}n, t_{j}+\f {\ga \be \da}n]} \(|\sin t|+\f
{4\be\ga \da} n\)
^p \\
&\    \     \  \hspace{10mm} \times \( \int_{t_{j}-\f {\ga \be
\da}n}^{t_{j}+\f {\ga \be \da}n} |T'(g(t))| ^p
 |\sin t| \, dt \) \( \int_{t_j-\f {\ga  \da}n}^{t_{j}+
 \f {\ga  \da}n} W_\al (t) \,   dt \)
 \\
&\leq C\be^{2p-1}  \al^{p+1} \(\f\da n\)^p  \int_{t_{j}-\f {\ga
\be \da}n}^{t_{j}+\f {\ga \be \da}n} | T'(g(t))|^p  |\sin t|
\(|\sin t|+\f 1 n\) ^p W_{ n_1,\al} (t)\, dt,\end{align*} where in
the first inequality we have used ($\ref{2-13-a}$),
($\ref{3-5-b}$) and H\"older's inequality, and in the last
inequality we have used the doubling property of $W_\al$. Thus, by
($\ref{3-4-b}$), we deduce
\begin{align*} \sum_{j=1}^{m_n}& \( \max_{x, y \in B_{\rho_1} ( \xi_j, \f{\be\da} n)
 }
|T(x)-T(y)|^p\) \(\int_{B_{\rho_1} ( \xi_j, \f\da n)} W(\xi)\, d\xi \)\\
& \leq C\be^{2p}  \al^{p+1} \(\f \da n\)^p \int_{0}^\pi  |
T'(g(t))|^p \sin t \(\sin t+\f 1 n\) ^p W_{ n_1,\al} (t)\,
dt\\
&\sim \be^{2p}\(\f \da n\)^p \int_{-\al}^\al  |T'(x)|^p  \(
\sqrt{\al^2-x^2} + \f \al n  \)^p  \widetilde{W}(x)\, dx \equiv
:I,\ \ \    \  \text{(by ($\ref{2-13-a}$))}
 \end{align*}
where $\widetilde{W}(x)= W_{ n_1,\al} (\arccos (\sin x
/\sin\al))$. Note that  $\wt{W}_\al(t)\equiv \widetilde{W}(g(t)) =
W_{ n_1,\al}(t)$ and that  $W_{ n_1,\al}(t)$ is a doubling weight
on $\mathbb {S}^1$ with the doubling constant depending only on
that of $W_\al$.   It follows by ($\ref{2-13-a}$) and \mbox{Lemma
3.1} that
\begin{align*}
I&\leq (C\be^2\da)^p \int_{-\al}^\al  |T(x)|^p \widetilde{W}(x)\,
dx \sim \be^{2p}\da^p \al \int_{0}^\pi  |T_\al (t)|^p   W_{
n_1,\al} (t) \sin
t \, dt\\
&\sim \be^{2p}\da^p \al\int_{0}^\pi  |T_\al (t)|^p   W_{\al} (t)
\sin t \, dt\sim \be^{2p}\da^p \int_{-\al}^\al  |T(x)|^p   W(x)\,
dx.\end{align*} This completes the proof of Lemma 3.2.
\end{proof}

\begin{lem}Let $W$
be a doubling weight on $\sph$. Let    $\da\in (0,1)$, $\be \ge 1$
and $0<p<\infty$.  Suppose that  $n$ is a positive integer and
$\Ld\subset \sph$ is $\f\da n$-separable with respect to the
geodesic metric $d(\cdot,\cdot)$ on $\sph$.  Then for all
$f\in\Pi_n^d$,
\begin{equation}\label{3-2} \sum_{\og\in\Ld} \(\max_{x\in B(\og,
\f {\be \da} n)} |f(x)|^p\)W\(B(\og, \da/n)\)
 \leq C \int_{\sph}
|f(x)|^p W(x)\, d\sa(x),\end{equation} and
\begin{equation}\label{3-3} \sum_{\og\in\Ld} \(\max_{x,y\in B(\og,
\f {\be\da} n)} |f(x)-f(y)|^p\)W\(B(\og, \da/n)\)
 \leq (C\da)^p \int_{\sph}
|f(x)|^p W(x)\, d\sa(x),\end{equation} where $C$ depends only on
$d$, $\be$,  $p$ and the doubling constant of $W$.\end{lem}
\begin{proof}  Equation ($\ref{3-2}$) is a direct consequence of Equation ($\ref{3-3}$).
Equation ($\ref{3-3}$) with $\be =1$ was proved in [D1, Corollary
3.3], and the proof there works equally well for  $\be >
1$.\end{proof}

Our fourth   lemma is due to Mhaskar, Narcowich and  Ward [MNW,
Proposition 4.1].   Let $X$ be a finite dimensional normed linear
space, $X^\ast$ be its dual, and $Z\subset X^\ast$ be a finite
subset with cardinality $m$.  We say $Z$ is a norming set for $X$
if the operator $x\mapsto (y^\ast (x))_{y^\ast\in Z}$ from $X$ to
$\mathbb{R}^m$ is injective. A functional $x^\ast\in X^\ast$ is
said to be positive with respect to $Z$ if for all $x\in X$,
 $x^\ast (x) \ge 0$ whenever $\dmin _{y^\ast \in Z}  y^\ast (x) \ge 0$.

\begin{lem}  {\rm ([MNW]). }\    \   Let $X$ be a
finite dimensional normed linear space, $X^\ast$ be its dual,
$Z\subset X^\ast$ be a finite, norming set for $X$, and $x^\ast
\in X^\ast$ be positive with respect to $Z$. Suppose further that
$\dsup_{x\in X} \dmin_{y^\ast \in Z} y^{\ast} (x)
>0$. Then there exists a sequence of  nonnegative numbers
$\ell_{y^\ast}$, ($y^{\ast}\in Z$) such that for any $x\in X$,
$$ x^\ast (x) =\sum_{y^\ast \in Z} \ell_{y^\ast} y^\ast (x).$$
\end{lem}

Recall that for $0<a<b\leq \pi$ and $e\in \sph$,
$$ B(e; a, b) =\Bl\{ y\in\sph:\   \  a\leq d(e, y) \leq b\Br\}.$$
Our final lemma, Lemma 3.5 below, will play a crucial role in the
proof of Theorem 1.1.

\begin{lem}
Let $\be \ge 1$, $\al\in (0,\f12]$ and $1\leq p<\infty$. Let $\Ld$
be a $(\f \da n, \rho)$-separable subset of $B(e,\al)$. Then for
all $f\in\Pi_n^d$,
\begin{equation*}
\sum_{\og\in\Ld\cap B(e; \f \al{12},\al)}\( \max_{x,y\in
B_\rho(\og, \f{\be \da}n)} |f(x)-f(y)|^p\) | B_\rho (\og, \f\da
n)|\leq  (C \da)^p \int_{B(e,\al)} |f(x)|^p \,
d\sa(x),\end{equation*} where $C>0$ depends only on $d$, $p$ and
$\be$.
\end{lem}

For the moment, we take Lemma 3.5 for granted and proceed with the
proof of our main results.

\vspace{5mm}

 {\it Proof of Theorem 1.1.}\   \      Let $1\leq p<\infty$, $\be \ge 1$, $f\in\Pi_n^d$ and let
$\Ld$ be a $(\f \da n, \rho)$-separable subset of  $B(e,\al)$. Set
$\Ld_1 =\Ld \bigcap B(e, \f \al {12})$. Then by Lemma 3.5, it will
suffice to show that
\begin{equation}\label{3-4}
\sum_{\og\in\Ld_1}\( \max_{x,y\in B_\rho(\og, \f{ \be \da}n)}
|f(x)-f(y)|^p\) | B_\rho (\og, \f\da n)|\leq  (C \da)^p
\int_{B(e,\al)} |f(x)|^p \, d\sa(x).\end{equation}

 For the proof of ($\ref{3-4}$), we take $u\in B(e, \al)$ so that $d(u, e) = \f
\al 6$.  Then associated with  the spherical cap $B(u, \al/2)$, we
define
$$ \wt{\rho}(x,y)\equiv \rho_{ _{B (u, \al/2)}}(x, y)  =
\f 2{\al}\sqrt{ (d(x, y))^2 +\f\al2 \Bl|
\sqrt{\wt{b}_x}-\sqrt{\wt{b}_y}\Br|^2},$$ where  $x, y\in
B(u,\f\al2)$,  and $\wt{b}_x\equiv b_{x, B(u, \al/2)}$ denotes the
shortest distance from $x\in B(u,\al/2)$ to the boundary of
$B(u,\f \al2)$.
   Since
$$\Ld_1\subset B(e, \f \al {12})\subset B(u; \f \al{12}, \f \al 4),$$
  by Lemma 2.2 applied  to both  $\rho\equiv
\rho_{_{B(e,\al)}}$ and $\wt{\rho}\equiv \rho_{_{B(u,\al/2)}}$, we
conclude that the following statements hold
  true:
\begin{align}
\wt{\rho} (x, y)&\sim \rho(x,y)\sim \f {d(x,y)}\al,\ \ \text{for
any $x, y\in\Ld_1$},\     \    \  \text{(by ($\ref{2-24-a}$))}\label{3-5}\\
|B_{\rho}(\og, r)|&\sim |B_{\wt{\rho}}(\og, r)|\sim (\al r)^{d}, \
\ \ \text{for any
$\og\in\Ld_1$ and $r\in (0,1)$}, \    \    \  \text{(by ($\ref{2-28-a}$))}\label{3-6}\\
B_{\rho}(\og, r)&\subset B_{\wt{\rho}}(\og, 2C_2r),\ \ \ \text{for
any $\og \in\Ld_1$ and $r>0$},\    \   \  \text{(by
($\ref{2-27-a}$))}\label{3-7}\end{align} where  all constants of
equivalence depend only on $d$, and $C_2$ is the absolute constant
in Equation ($\ref{2-27-a}$) with $\va=\f12$.   By Equation
($\ref{3-5}$), we know that there exists an absolute constant
$\ga\in (0,1)$ such that $\Ld_1$ is
 $(\f {\ga\da} { n}, \wt{\rho})$-separable  in $B(u, \al/2)$.   However, on the other hand,    using
($\ref{3-6}$) and ($\ref{3-7}$), we deduce  that the sum on the
left-hand side of ($\ref{3-4}$) is controlled by
$$J\equiv C \sum_{\og\in \Ld_1 \cap B(u;  \f \al {24},  \f \al2)}\(
\max_{x,y\in B_{\wt{\rho}}(\og, \f{A \ga \da}n)} |f(x)-f(y)|^p\) |
B_{\wt{\rho}} (\og, \f{\ga \da} n)|, $$ where $A=2\be C_2/\ga$.
Therefore, by  Lemma 3.5  applied to  $B(u,\f \al 2)$ and
$\wt{\rho}\equiv \rho_{_{B(u, \al/2)}}$,  it follows that
$$ J\leq C (\ga \da)^p \int_{B(u,\al/2)} |f(x)|^p \, d\sa(x)\leq (C\da)^p
\int_{B(e,\al)} |f(x)|^p \, d\sa(x),$$ which proves  the desired
equation ($\ref{3-4}$) and hence ($\ref{1-6}$).\hb

\vspace{5mm}

Now we turn to the proofs of Corollaries. \vspace{5mm}

{\it Proof of Corollary 1.2.}\     \  Let $C_1$ denote  the
constant in ($\ref{1-6}$) with $\be=p=1$. Let  $\da \in (0,\f
1{4C_1}]$, $n_1=[n/(4C_1 \da)]$ and let $\Ld$ be a maximal $(\f\da
n, \rho)$-separable subset of $B(e,\al)$.   We shall prove  that
there exists a sequence of positive numbers $\ld_\og$, $\og\in\Ld$
such that $\ld_\og \sim |B_{\rho} (\og, \f \da n)|$ and
($\ref{1-1}$) holds for all $f\in \Pi_{n_1}^d$.
 The  idea of our proof below  is from [NPW].

 Note,  by Lemma 2.2 (v), that
 \begin{equation}\label{3-8}
1\leq B(x)\equiv \sum_{\og\in \Ld} \chi_{_{B_{\rho} (\og, \f\da
n)}}(x) \leq C,\    \    \    \  x\in B(e,\al),\end{equation}
where $C\ge 1$ depends only on $d$. We  define  the following
linear functional on $\Pi_{n_1}^d$:
$$\ell(f) =2\int_{B(e,\al)} f(x)\, d\sa(x)  -\sum_{\og\in\Ld}
 \(\int_{B_{\rho} ( \og,  \f\da n)}\f{d\sa(x)}{B(x)} \) f(\og),\    \    \
 f\in\Pi^d_{n_1}.$$We then claim
that there exists a sequence of nonnegative numbers $\mu_\og$,
$\og\in\Ld$ such that
\begin{equation}\label{3-9}
\ell(f) =\sum_{\og\in\Ld} \mu_\og f(\og),\     \      \ \text{for
all $f\in \Pi_{n_1}^d$}.\end{equation}

  For the proof of the claim ($\ref{3-9}$), we note  that, by ($\ref{1-6}$),  each
$f\in\Pi^d_{n_1}$ is uniquely determined by its restriction to the
set $\Ld$. (This can  be also seen from the proof below.) Thus, in
view of Lemma 3.4, it will suffice to prove that for any
$f\in\Pi_{n_1}^d$ with $\dmin_{\og\in\Ld} f(\og) \ge 0$, $\ell(f)
\ge 0.$ To see this, we note that if $f(\og)\ge 0$, then for any $
x\in B_{\rho} (\og, \f \da n)$,
\begin{align*}
2f(x)-f(\og)&\ge \max_{z\in B_{\rho} (\og, \f \da n) } |f(z)|
-\max _{ z\in B_{\rho} (\og, \f \da n)} \Bl( |f(z)|-f(\og)
-2f(x)+2f(\og) \)\\
&\ge \max_{z\in B_{\rho} (\og, \f \da n) }  |f(z)| - 3\max_{y\in
B_{\rho} (\og, \f \da n) }
  |f(y)-f(\og)|.\end{align*}
 Thus,  for  $f\in \Pi_{n_1}^d$
with $\dmin_{\og\in \Ld} f(\og) \ge 0$, we have
\begin{align*}
&\ell(f) =\sum_{\og\in\Ld} \int_{B_{\rho}
(\og,  \f \da n)}(2f(x) -f(\og)) \f{ d\sa(x)}{B(x)}\notag\\
& \ge \sum_{\og\in\Ld}  \Bl[ \max_{z\in B_{\rho} (\og, \f \da n) }
|f(z)| - 3\max_{y\in B_{\rho} (\og, \f \da n) }
  |f(y)-f(\og)|\Br] \int_ {B_{\rho}
(\og, \f \da n)}\f{ d\sa(x)}{B(x)} \notag\\
&\ge  \int_{B(e,\al)} |f(x)|\, d\sa(x)-  3\sum_{\og\in\Ld}
\(\max_{y \in B_{\rho} (\og, \f \da n) }
  |f(y)-f(\og)|\) |B_{\rho}
(\og, \f \da n)|,
\end{align*}
which,  by ($\ref{1-6}$) with $n$  replaced by $n_1$,  is greater
 or equal
\begin{equation*}
(1- 3C_1 \f {n_1 \da } {n})\int_{B(e,\al)}|f(x)|\, d\sa(x)\ge \f
14 \int_{B(e,\al)} |f(x)|\, d\sa(x) \ge 0.\end{equation*} This
proves the claim ($\ref{3-9}$).

Now setting
$$\ld_\og=\f12\mu_\og+ \f12\int_{B_{\rho} (\og, \da/n)}
\f{d\sa(x)}{B(x)},\     \    \    \og\in\Ld,$$  and taking into
account ($\ref{3-8}$) and ($\ref{3-9}$),
 we conclude that   ($\ref{1-1}$)  with
$\ld_\og$ satisfying the condition $\ld_\og \ge C^{-1} |B_{\rho}
(\og, \f \da n)|$  holds for all $f\in \Pi_{n_1}$. Thus, it
remains to show the inequality
\begin{equation}\label{3-13-b}
\ld_\og \leq  C |B_{\rho} (\og, \f \da n)|,\   \
\og\in\Ld,\end{equation} where $C>0$ depends only on $d$. To this
end,   we set   $ n_2 =[ n_1/(2d+2)]$ and
$$ A_{n_1} (\cos t) = \ga_{n_2} \( \f { \sin ( n_2+\f12) t}
{\sin\f t2}\)^{2d+2},\   \   t\in[-\pi,\pi],$$ where $\ga_{n_2}$
is a positive constant chosen so that $A_{n_1} (1)=1$.   Then it
is easy to verify  that
\begin{equation}\label{3-11}
|A_{n_1} (\cos t) | \leq C ( 1+ n_1 |t|) ^{-2d-2} ,\    \     \
t\in [-\pi, \pi].\end{equation}

 Now for a fixed  $\og =(
\og_1,\cdots,\og_{d},\og_{d+1})\equiv (\og', \og_{d+1})\in\Ld$, we
define
\begin{align*}
 f_{n_1} ( y)& = A_{n_1}\( \f {y'\cdot \og'} {\sin^2\al}  + \f {
  \sqrt{ y_{d+1}^2-\cos^2\al}\sqrt{
\og_{d+1}^2-\cos^2\al}} {\sin^2\al} \)\\
&\    \    \    \  + A_{n_1}\( \f {y'\cdot \og'} {\sin^2\al}  -\f
{  \sqrt{ y_{d+1}^2-\cos^2\al}\sqrt{ \og_{d+1}^2-\cos^2\al}}
{\sin^2\al} \),\end{align*} where $y=(y_1,\cdots, y_{d},
y_{d+1})\equiv ( y',y_{d+1})\in B(e,\al)$. Since $A_{n_1}$ is an
algebraic polynomial of degree at most $n_1$ on $[-1,1]$, it
follows that $f_{n_1} \in\Pi_{n_1}^d$. Note, on the other hand,
\begin{align*}
 \arccos  &\Bl[ \f {y'\cdot \og'} {\sin^2\al}  \pm \f {
 \sqrt{ y_{d+1}^2-\cos^2\al}\sqrt{ \og_{d+1}^2-\cos^2\al}} {\sin^2\al}
\Br]\\[3mm]
&\sim \f 1{\sin\al} \sqrt{ |y'-\og'|^2+\bl|
\sqrt{y_{d+1}^2-\cos^2\al}\mp \sqrt{\og_{d+1}^2-\cos^2\al}\br|^2}\\[3mm]
 &\ge  \rho_5(y,\og)\ge C \rho(y,\og),\end{align*}
 where the last two inequalities follow  by ($\ref{2-21-a}$) and
 ($\ref{2-22-a}$), respectively.
Thus, by ($\ref{3-11}$), we obtain
\begin{equation*}
0\leq f_{n_1}( y)\leq C ( 1+ n_1\rho(y,\og))^{-2d-2},\   \   \ \
 y\in B(e,\al).\end{equation*} Now applying the cubature formula
 ($\ref{1-1}$) to $f_{n_1}$, we deduce
\begin{align*}
\ld_\og& \leq \ld_\og f_{n_1}(\og) \leq \sum_{\xi\in \Ld} \ld_\xi
f_{n_1}(\xi)=\int_{B(e,\al)} f_{n_1}(y)\, d\sa(y)\\
& \leq C \sum_{j=0}^{\infty} \int_{\{y\in B(e,\al): \    \  \f
j{n_1} \leq \rho(\og, y)\leq \f{j+1}{n_1}\}} f_{n_1}(y)\,
d\sa(y)\\
&\leq C \bl|B_{\rho} (\og, \f 1{n_1})\br|\sum_{j=0}^\infty
(j+1)^{-d-1}\leq C |B_{\rho} (\og, \f \da n)|,\end{align*} which
gives ($\ref{3-13-b}$) and hence completes the proof of Corollary
1.2.\hb\\

{\it Proof of Corollary 1.3.}\    \   Let $\mathcal{A}$ be a
maximal $(\f 1n,\rho)$-separable subset of $B(e,\al)$.  Then by
Equation ($\ref{1-6}$) and Lemma 2.2 (v), it is easily seen   that
for $f\in \Pi_n^d$ and $1\leq p<\infty$,
\begin{equation}\label{3-15-Feb}\sum_{\xi\in \mathcal{A}} \( \max_{ x\in B_{\rho} (\xi,
\f 1n)} |f(x)|^p \) \Bl|B_\rho (\xi, \f 1n)\Br|\leq C
\int_{B(e,\al)} |f(x)|^p \, d\sa(x).\end{equation} Using this last
fact, we obtain
\begin{align*}
&\sum_{\og\in\Ld} |f(\og)|^p \Bl| B_\rho (\og, \f 1n)\Br| \leq
\sum_{\xi\in \mathcal{A}} \sum_{\og \in\Ld\cap  B_{\rho}(\xi, \f
1n)} |f(\og)|^p
\Bl|B_\rho ( \og, \f 1n)\Br|\\
&\leq C\sum_{\xi\in\mathcal{A}} \(\max_{x\in B_\rho(\xi, \f 1n)}
|f(x)|^p
\)\# \( \Ld \cap B_\rho(\xi, \f 1n)\)\Bl|B_\rho ( \xi, \f 1n)\Br|\\
& \leq C \tau \int_{B(e,\al)} |f(x)|^p \, d\sa(x),\end{align*}
where in the second inequality we have used the fact that
$\Bl|B_\rho ( \xi, \f 1n)\Br|\sim \Bl|B_\rho (x, \f 1n)\Br|$
whenever $x\in B_\rho ( \xi, \f 1n)$, and in the last  inequality
we have used ($\ref{3-15-Feb}$) and the definition of $\tau$.
This completes the proof of Corollary 1.3.

\vspace{5mm}

{\it Proof of Corollary 1.4.}\   \    Corollary 1.4 is a simple
consequence of Equation ($\ref{1-6}$) and Lemma 2.2 (v). We omit
the detail. \hb

\vspace{5mm}

Now it remains to show Lemma  3.5.

\vspace{5mm}

{\it Proof of Lemma 3.5.}\   \    Suppose that $f\in\Pi_n^d$ and
$\Ld$ is a $(\f \da n, \rho)$-separable subset of $B(e,\al)$.  We
set $\Ld_2=\Ld\cap B(e; \f \al{12},\al)$. Since Lemma 3.5 is a
direct consequence of \mbox{Lemma 3.2} in the case when $d=1$,  we
shall assume $d\ge 2$ in the proof below.  Also, without loss of
generality we may assume that $e=(0,0,\cdots,0, 1)\in\sph$.

Recall that $\rho_4$ is a metric on $B(e,\al)$ defined by
($\ref{2-20-a}$). It follows by ($\ref{2-23-a}$) and
($\ref{2-26-a}$) with $\va =\f 1{24}$ that there exists an
absolute constant $C_3\ge 1$ such that
\begin{equation}\label{3-12} C_3^{-1} \rho (x, y) \leq  \rho_4 (x, y)\leq C_3 \rho (x,
y),\    \    \   \  \text{for all $x, y\in B(e;\f
\al{24},\al)$},\end{equation} and
\begin{equation}\label{3-13}
B_{\rho_4} (x, C_3^{-1} r)\subset B_{\rho} (x,  r)\subset
B_{\rho_4} (x, C_3 r),\    \  \text{for all $x\in B(e;\f \al{24},
\al)$ and $r>0$.}
\end{equation}

Next, recall that  $\rho_1\equiv \rho_{_{[-\al,\al]}}$  is  the
metric on $[-\al, \al]$ defined by ($\ref{2-1-a}$). Let
$\{v_i\}_{i=0}^{L_n}$ be a sequence of numbers in $[\f \al {12},
\al]$ satisfying the conditions  $\dmin _{0\leq i\neq j \leq L_n}
\rho_1 ( v_i, v_j)\ge \f {3^{-1}C_3^{-1}\da}n$ and $\displaystyle
[\f \al {12}, \al] \subset \bigcup_{i=0}^{L_n} B_{\rho_1} ( v_i,
\f {3^{-1}C_3^{-1}\da}n)$.
  Let $\{\xi_j\}_{j=0}^{M_n}$ be
a  maximal $(\f {3^{-1}C_3^{-1}\da}n,
d_{_{\mathbb{S}^{d-1}}})$-separable subset of  $\mathbb{S}^{d-1}$,
where $d_{_{\mathbb{S}^{d-1}}}$ denotes the usual  geodesic metric
on $\mathbb{S}^{d-1}$.  Set
$$\og_{ij}=(\xi_j\sin v_i,\cos v_i),\   \
0\leq i\leq L_n,\   \    \  0\leq j \leq M_n.$$ Then it's easily
seen that \begin{equation}\label{3-14} B(e; \f \al {12},\al)
\subset \bigcup_{i =0}^{L_n}\bigcup_{j=0}^{M_n} B_{\rho_4} (
\og_{ij}, \f {3^{-1}C_3^{-1}\da } n).\end{equation}
 On the other hand, by ($\ref{3-12}$) it follows that  $\Ld_2$
 is $(\f {C_3^{-1} \da} n,  \rho_4)$-separable. This means that
   \begin{equation}\label{3-15} \#  \(\Ld_2 \bigcap B_{\rho_4} ( \og_{ij}, \f{3^{-1} C_3^{-1} \da}n)\)\leq
 1,\     \   \ 0\leq i\leq L_n,\   \    \  0\leq j \leq
 M_n.\end{equation}
  Also, note that  if $\og \in \Ld_2 \bigcap B_{\rho_4} ( \og_{ij},
 \f{3^{-1} C_3^{-1}\da}n)$, then by ($\ref{3-13}$), ($\ref{2-28-a}$) and ($\ref{2-29-a}$),  for any given
 $\be \ge 1$,
 \begin{equation}\label{3-16}B_{\rho} (\og, \f{\be  \da }
 n)\subset B_{\rho_4} (\og_{ij}, \f{ C_4\da } n),\   \   \text{and}\   \
 | B_\rho (\og, \f\da n)|\sim | B_{\rho_4} (\og_{ij}, \f\da
 n)|,\end{equation}
 where  $C_4 =C_3\be + 3^{-1} C_3^{-1}.$

Therefore, setting  $$ \mathcal{A}=\Bl\{(i, j)\in [0, L_n]\times
[0, M_n]: \   \    \Ld_2 \bigcap B_{\rho_4} ( \og_{ij},
\f{3^{-1}C_3^{-1}\da}n)\neq \emptyset\Br\}, $$ and taking into
account ($\ref{3-14}$) and ($\ref{3-15}$), we conclude that for
every $\og\in\Ld_2$, there exists a unique $(i, j)\in\mathcal {A}$
for which  $\og \in \Ld_2 \bigcap B_{\rho_4} ( \og_{ij},
 \f{3^{-1} C_3^{-1}\da}n)$ and  ($\ref{3-16}$) holds. This implies
\begin{align}&\sum_{\og\in\Ld_2}\( \max_{x,y\in B_\rho(\og,
\f{\be\da}n)} |f(x)-f(y)|^p\) | B_\rho (\og, \f\da
n)|\notag\\
&\leq  C  \sum_{(i,j)\in\mathcal{A}}  \( \max_{x\in B_{\rho_4} (
\og_{ij}, \f {C_4\da}n )} |f(x)-f(\og_{ij})|^p \) |B_{\rho_4}
(\og_{ij}, \f \da n)|\notag\\
&\leq C\sum_{i=0}^{L_n}\sum_{j=0}^{M_n}  \( \max_{x\in B_{\rho_4}
( \og_{ij}, \f {C_4\da}n )} |f(x)-f(\og_{ij})|^p \) |B_{\rho_4}
(\og_{ij}, \f \da n)|\equiv \Sigma. \label{3-17}\end{align} Thus,
the proof of Lemma 3.5  is now reduced to the proof of the
following inequality:
\begin{equation}\label{3-18}
\Sigma\leq (C\da)^p \int_{B(e,\al)} |f(x)|^p\,
d\sa(x),\end{equation} where $\Sigma$ is defined by
($\ref{3-17}$), and the constant $C$ depends only on $d$, $p$ and
$\be$.

 For  the rest of the proof,
we shall write $\dsum_{i,j}$ for
$\dsum_{i=0}^{L_n}\dsum_{j=0}^{M_n}$, $\dsum_i$ for
$\dsum_{i=0}^{L_n}$, and  $\dsum_j$ for $\dsum_{j=0}^{M_n}$.
Moreover, given  $r>0$ and  $\xi\in \mathbb{S}^{d-1}$,  we denote
by  $ B(\xi, r)\equiv B_{_{d_{\mathbb{S}^{d-1}}}}(\xi,r)$  the
spherical cap $\{\eta\in \mathbb{S}^{d-1}: \ \ \arccos \xi\cdot
\eta \leq r\}$ in $\mathbb{S}^{d-1}$.

  To
show ($\ref{3-18}$), we define $ g(v, \eta) = f( \eta\sin v, \cos
v),$ where  $\eta\in\mathbb{S}^{d-1}$ and $v\in[-\al,\al]$, and we
let $F$ be a polynomial on $\mathbb{R}^{d+1}$ of total degree at
most $n$ whose restriction to $\sph$ is $f$. Then, by the chain
rule,  we have,  for $\eta=(\eta_1,\cdots,\eta_{d})\in
\mathbb{S}^{d-1}$,
$$ \f {\p g(v, \eta)}{\p v} =\sum_{k=1}^{d} \f{\p F(\eta \sin v,
\cos v)}{\p x_k} \eta_k \cos v-\f{\p F(\eta \sin v, \cos v)}{\p
x_{d+1}} \sin v.$$   It follows  that  $\df {\p g(\cdot, \eta)}{\p
v} \in\Pi^1_n$ for each  fixed $\eta \in\mathbb{S}^{d-1}$, and
$\df {\p g(v, \cdot)}{\p v} \in \Pi_n^{d-1}$  for each  fixed
$v\in [-\al,\al]$.

   Now for each $(i, j)\in [0, L_n]\times [0, M_n]$,  we assume
$$ \max_{x\in B_{\rho_4} ( \og_{ij}, \f {C_4\da}n )}
|f(x)-f(\og_{ij})|=|f(x^\ast_{ij})-f(\og_{ij})|,$$ where
$x_{ij}^\ast=( \xi^\ast_{ij} \sin \ta^\ast _{ij},\cos
\ta^\ast_{ij})\in B_{\rho_4} ( \og_{ij}, \f {C_4\da}n ) $, that
is,   $\xi^\ast_{ij} \in B(\xi_j, \f {C_4\da}n)$,
$\ta_{ij}^\ast\in B_{\rho_1} (v_i, \f {C_4\da}n)$. Then we have
\begin{align*}
&|f(x^\ast_{ij})-f(\og_{ij})|^p= \bl|g(\ta^\ast_{ij},
\xi^\ast_{ij})
-g(v_i, \xi_j)\br|^p\\[4mm] &\leq  2^p
\bl|g(v_i, \xi_j)-g(v_i, \xi^\ast_{ij})\br|^p+2^p \bl|g(v_i,
\xi^\ast_{ij})-g(\ta^\ast_{ij},
\xi^\ast_{ij})\br|^p\\[4mm]
&\leq 2^p  \max_{\eta\in B(\xi_j, \f {C_4\da} n) } | g(v_i,
\eta)-g(v_i, \xi_j)|^p+ 2^p \Bl|B_{\rho_1} (v_i, \f
{C_4\da}n)\Br|^{p-1}\int_{B_{\rho_1} (v_i, \f {C_4\da}n)}
\Bl|\f {\p g(v, \xi^\ast_{ij})} { \p v}\Br|^p \, dv \\
&\leq   2^p  \max_{\eta\in B(\xi_j, \f {C_4\da} n) } | g(v_i,
\eta)-g(v_i, \xi_j)|^p    +    C \bl|B_{\rho_1} (v_i, \f
{\da}n)\br|^{p-1}\int_{B_{\rho_1} (v_i, \f {C_4\da}n)} \Bl|\f {\p
g(v,
\xi_j)} { \p v}\Br|^p \, dv \\
&\   \    \     \   \   +       C    \bl|B_{\rho_1} (v_i, \f
{\da}n)\br|^{p-1}\int_{B_{\rho_1} (v_i, \f {C_4\da}n)}
\(\max_{\eta\in B(\xi_j, \f {C_4\da} n) } \Bl|\f {\p g(v, \xi_j)}
{ \p v}-\f {\p g(v, \eta)} { \p v}\Br|^p\)
\, dv\\
 &\equiv A_{ij}+B_{ij}+C_{ij},\end{align*}
 where  in the second inequality we have used H\"older's
 inequality and the fact that $B_{\rho_1}(v_i, \f {C_4\da}n)$ is a
 subinterval of $[-\al,\al]$ containing  both $v_i$ and
 $\ta_{ij}^\ast$, and in the third inequality we have used
 ($\ref{2-6}$). Since, for any $\da' \in (0,\da)$,
 $\Ld_2$ is, again,  $(\rho,\f{\da'}n)$-separable,  without loss of generality, we may assume $\da \in
 (0, \f 1{ 24 C_4})$.
 Thus,  by Lemma 2.1 (ii), we deduce
   \begin{equation}\label{3-22-a}B_{\rho_1} ( v, \f {C_4\da} n) \subset
   [ v- \f { C_4\da
 \al}{n}, \al]\subset [ \f \al{24}, \al],\    \    \  \text{for any
 $v\in [\f \al{12},\al]$}.\end{equation}  Hence, for each  $(i,j)$,
 we have
\begin{equation*}
\bl|B_{\rho_4}(\og_{ij},\f \da n)\br|= C_d\bl|B(\xi_j, \f \da
n)\br| \int_{B_{\rho_1} (v_i, \f \da n)} \sin^{d-1}\ta\, d\ta\sim
\(\f{\da\al} n\)^{d-1} |B_{\rho_1} (v_i, \f \da
n)|.\end{equation*}
 It follows by ($\ref{3-17}$) that
\begin{align}
\Sigma &\leq C\(\f{\da\al} n\)^{d-1}\sum_{i, j} A_{ij} |B_{\rho_1}
(v_i, \f \da n)|+C\(\f{\da\al} n\)^{d-1}\sum_{i, j}
B_{ij}|B_{\rho_1} (v_i, \f \da n)|\notag \\
&  \    \    \     \   +C\(\f{\da\al} n\)^{d-1}\sum_{i, j} C_{ij}
|B_{\rho_1} (v_i, \f \da n)|
 \notag\\
&\equiv \Sigma_1+\Sigma_2+\Sigma_3.\label{3-24-a}\end{align}

 For the first sum $\Sigma_1$, we have
\begin{align*} \Sigma_1&\leq C\al^{d-1} \sum_i
|B_{\rho_1} (v_i, \f {\da}n)| \Bl[ \(\f \da n\)^{d-1}\sum_j
\max_{\eta\in B(\xi_j, \f {C_4\da} n) } | g(v_i,
\eta)-g(v_i, \xi_j)|^p\Br]  \\
&\leq (C \da)^p \sum_i \al^{d-1} |B_{\rho_1} (v_i, \f {\da}n)|
\int_{\mathbb{S}^{d-1}}
|g(v_i, \eta)|^p\, d\sa(\eta) \\
&\leq (C\da)^p \int_{\mathbb{S}^{d-1}}\Bl[ \sum_i |g(v_i,\eta) |^p
\int_{B_{\rho_1} (v_i,  \f \da n)} |\sin^{d-1} v|\, dv\Br] \, d\sa(\eta) \\
&\leq (C\da)^p \int_{\mathbb{S}^{d-1}} \int_{-\al}^\al
|g(v,\eta)|^p |\sin^{d-1} v| \, dv\,
d\sa(\eta)=(C\da)^p\int_{B(e,\al)} |f(x)|^p\, d\sa(x),\end{align*}
where in the second inequality we have used Lemma 3.3
($\ref{3-3}$) and the fact that $g(v_i, \cdot)\in\Pi_n^{d-1}$ for
each fixed $i$, in the third inequality we have used
($\ref{3-22-a}$),  and the last inequality follows by Lemma 3.2
and Lemma 2.1 (iv).

For the second sum $\Sigma_2$, we have
\begin{align*}
\Sigma_2&\leq C \al^{d-1} \sum_i \Bl|B_{\rho_1} (v_i, \f
{\da}n)\Br|^{p}\int_{B_{\rho_1} (v_i, \f {C_4\da}n)} \Bl[\(\f\da
n\)^{d-1}\sum_j \Bl|\f {\p g(v, \xi_j)} { \p v}\Br|^p\Br] \, dv\\
&\leq C  \sum_i \int_{B_{\rho_1} (v_i, \f {C_4\da}n)}
|\sin^{d-1}v| \Bl|B_{\rho_1} (v, \f {\da}n)\Br|^{p}
\(\int_{\mathbb{S}^{d-1}}\Bl| \f {\p
g(v, \xi)} { \p v}\Br|^p\, d\sa(\xi)\)\, dv\\
&\leq C\(\f \da n \)^p \int_{-\al}^\al |\sin^{d-1} v| \( \f{ \al}n
+\sqrt{\al^2-v^2} \)^p \(\int_{\mathbb{S}^{d-1}}
\Bl| \f {\p g(v, \xi)} { \p v}\Br|^p \, d\sa(\xi)\)\, dv\\
&=C\(\f \da n \)^p\int_{\mathbb{S}^{d-1}} \(\int_{-\al}^\al
|\sin^{d-1} v| \( \f{ \al}n +\sqrt{\al^2-v^2} \)^p\Bl| \f {\p
g(v, \xi)} { \p v}\Br|^p\, dv\)\, d\sa(\xi)\\
 &\leq (C\da)^p \int_{\mathbb{S}^{d-1}}  \Bl[ \int_{-\al}^\al |
g(v, \xi)|^p |\sin^{d-1} v|  \, dv \Br] \, d\sa(\xi)=(C\da)^p
\int_{B(e,\al)} |f(x)|^p\, d\sa(x),
\end{align*}
where in the second inequality, we have used ($\ref{3-22-a}$),
($\ref{2-6-a}$), ($\ref{3-2}$) and the fact that $\df {\p g(v,
\cdot)}{\p v}\in \Pi_n^{d-1}$ for each  fixed $v$, in the third
inequality, we have used ($\ref{2-7-a}$) and ($\ref{2-6}$), and in
the last inequality, we have used ($\ref{3-1}$) and the fact that
$ g(\cdot, \xi)\in \Pi_n^{1}$ for each fixed
$\xi\in\mathbb{S}^{d-1}$.

 For the third sum $\Sigma_3$, we have
\begin{align*}
\Sigma_3 &\leq C \al^{d-1} \sum_i |B_{\rho_1} (v_i, \f
{\da}n)|^{p}\\
&\    \   \times \int_{B_{\rho_1} (v_i, \f {C_4\da}n)} \Bl[\(\f
\da n\)^{d-1} \sum_j \(\max_{\eta\in B(\xi_j, \f {C_4\da} n) }
\Bl|\f {\p g(v, \xi_j)} { \p v}-\f {\p g(v,
\eta)} { \p v}\Br|^p\) \Br]\, dv\\
&\leq (C \da)^p \(\f \da n \)^p \sum_i \int_{B_{\rho_1} (v_i, \f
{C_4\da}n)} |\sin^{d-1} v | \( \f { \al}n
+\sqrt{\al^2-v^2}\)^p\\
&\    \    \  \times  \(\int_{\mathbb{S}^{d-1}} \Bl|\f {\p g(v,
\xi)} { \p v}\Br|^p \,
d\sa(\xi) \)\, dv\\
 &\leq C\da^p \(\f \da n \)^p \int_{\mathbb{S}^{d-1}}
\( \int_{-\al}^\al \Bl| \f {\p g(v, \xi)} { \p v}\Br|^p
|\sin^{d-1} v| \( \f{
\al}n +\sqrt{\al^2-v^2} \)^p \, dv\) \, d\sa(\xi)\\
&\leq C\da^{2p} \int_{\mathbb{S}^{d-1}}  \( \int_{-\al}^\al | g(v,
\xi)|^p |\sin^{d-1} v| \, dv \)\, d\sa(\xi)=C\da^{2p}
\int_{B(e,\al)} |f(x)|^p\, d\sa(x),
\end{align*}
where in the second inequality we have used ($\ref{3-22-a}$),
($\ref{2-6-a}$), ($\ref{2-6}$) ($\ref{3-3}$) and the fact that
$\df {\p g(v, \cdot)}{\p v}\in \Pi_n^{d-1}$ for a fixed $v$, in
the third inequality we have used Lemma 2.1 (iv),  and in the last
inequality, we have used ($\ref{3-1}$) and the fact that $g(\cdot,
\xi)\in \Pi_n^{1}$ for each fixed $\xi\in\mathbb{S}^{d-1}$.

Now putting the above estimates together, and taking into account
($\ref{3-24-a}$), we deduce the desired inequality ($\ref{3-18}$),
and hence complete the proof of Lemma 3.5.\hb

\section{Proofs of the main results  for $\al \in [\f 12, \pi)$}

  Let $\va\in (0,1)$. In this section we shall prove  Theorem 1.1 and Corollaries
  1.2--1.4
  in the case when $\al \in [\f 12,
  \pi-\va]$. It turns out that the main results in this case can be deduced
    from the already proven case $\al\in (0,\f 12]$.
 Without loss of generality we may assume  in this
section that $d\ge 2$ and  $e=(0,\cdots,0, 1)\in\sph$. (The proof
for the case $d=1$ is similar and in fact, much simpler.)     For
$x= (\eta\sin\ta, \cos\ta)$ with $\ta\in [0,\pi]$ and
$\eta\in\mathbb{S}^{d-1}$, we define
\begin{equation}\label{4-1}Tx: = \bl(\eta \sin (8\ta), \cos (8\ta)\br).\end{equation}Then $T$ is a map from
$B(e,\f\al 8)$ to $B(e,\al)$.   Also, we set
\begin{equation}\label{4-2} D(\cos\ta) :=\f {
\sin^{d-1}(8\ta)}{\sin^{d-1}\ta},\ \ \ta\in [0,\pi].\end{equation}
Then $D$ is an algebraic polynomial on $[-1,1]$ of degree
$7(d-1)$.

We need two lemmas, the first of which can be stated as follows.

\begin{lem}
Let $\al\in (0,\pi)$ and let  $T$  be defined by ($\ref{4-1}$).
Then the following statements hold true:

(i)\    \  If $f\in\Pi_{n}^d$ then $f{\small \circ} T\in
\Pi_{8n}^d$.

(ii)\     \  If $f$ is an integrable function on $B(e,\al)$, and
$D$ is the polynomial  defined by ($\ref{4-2}$), then we have
\begin{equation}\label{4-3}
\int_{B(e,\al)} f(x)\, d\sa(x) = 8\int_{B(e, \f \al 8)} f(Tx)
D(x\cdot e)\, d\sa(x).\end{equation} \end{lem}
\begin{proof}  We start with the proof of (i).  Setting
$$A(\cos\ta):=\cos 8\ta\    \   \text{and}\   \   \
B(\cos\ta):=\df{\sin (8\ta)}{\sin\ta},$$
 we obtain that  for $x=(\eta\sin\ta,\cos\ta)\equiv (x', x_{d+1})\in\sph$,
$$
f(Tx)=f(\eta\sin(8\ta),\cos(8\ta))= f\bl((\eta \sin\ta)
B(\cos\ta), A(\cos \ta)\br)=f( x' B(x_{d+1}), A(x_{d+1})).$$ Note,
however, that  $A$ is a polynomial on $[-1,1]$ of degree $8$, and
$B$ is a polynomial on $[-1,1]$ of degree $7$. Assertion (i) then
follows.

Next, we show (ii).  In fact, we have
\begin{align*}
\int_{B(e,\al)} f(x)\, d\sa(x) &=C_d \int_0^\al
\int_{\mathbb{S}^{d-1}} f(\eta \sin \ta, \cos\ta)\, d\sa(\eta)
\sin^{d-1} \ta \, d\ta \\
&= 8 C_d \int_0^{\al/8} \int_{\mathbb{S}^{d-1}} f(\eta \sin
(8\ta), \cos(8\ta))\, d\sa(\eta) \sin^{d-1} (8\ta) \, d\ta\\
&=8 C_d \int_0^{\al/8} \int_{\mathbb{S}^{d-1}} (f\circ T) (\eta
\sin
\ta, \cos\ta)\, d\sa(\eta) D(\cos\ta) \sin^{d-1} \ta \, d\ta\\
&=8\int_{B(e,\al/8)} f(Tx) D(x\cdot e) \, d\sa(x),\end{align*}
proving ($\ref{4-3}$).
\end{proof}

Let $T$ be the map from $B(e, \al/8)$ to $B(e, \al)$ defined by
($\ref{4-1}$) and let $T^{-1}$ denote its inverse.  Given a subset
$E$ of $B(e,\al)$, we write
$$ T^{-1} (E) =\Bl\{ x\in B(e, \al/8):\   \   Tx \in E\Br\}.$$
Also, we recall that $\rho_{_{B(e,\al)}}$ denotes the metric on
$B(e,\al)$ defined by ($\ref{1-3}$).  For simplicity, we shall
write $\rho_\al =\rho_{_{B(e,\al)}}$ and $
\rho_{\al/8}=\rho_{_{B(e,\al/8)}}$.

Now our second lemma can be stated as follows.

\begin{lem}
Let $\va\in (0,1)$ and $\al\in (0,\pi-\va]$.  Then there exists a
positive constant $C_5$ depending only on $d$ and $\va$ when $\va$
is small such that the following statements hold true:

(i)\   \  For any $x, y\in B(e, \al/8)$,
$$
C_5^{-1} \rho_{ \al/8} (x, y) \leq \rho_{\al} (Tx, Ty)\leq C_5
\rho_{ \al/8} (x, y).$$

(ii)\   \   For any $x\in B(e,\al)$ and $r>0$,
$$ B_{\rho_{\al/8}}(T^{-1} x, C_{5}^{-1}r) \subset T^{-1} \(
 B_{\rho_{\al}}( x, r)\) \subset B_{\rho_{\al/8}}(T^{-1} x,
 C_{5}r).$$

(iii)\   \ For any measurable subset $E$ of $B(e,\al)$,
$$ C_5^{-1} \bl|T^{-1} (E)\br|\leq |E| \leq C_5 \bl|T^{-1}
(E)\br|.$$

(iv)\   \  For any  $x\in B(e,\al)$ and $r\in (0, 1)$,
$$ C_5^{-1} \tr_{_{r, B(e,\al)}} (x) \leq |B_{_{\rho_\al}}(x, r)|\leq C_5
 \tr_{_{r, B(e,\al)}}
(x),$$ where $\tr_{_{r, B(e,\al)}}(x)$ is defined by
($\ref{1-5}$).
 \end{lem}
\begin{proof} (i)\   \   Let $x= (\eta\sin\ta, \cos \ta)$ and
$y=(\xi\sin t, \cos t)$ with $\ta, t\in [0,\al/8]$ and
$\xi,\eta\in \mathbb{S}^{d-1}$.  Then, by ($\ref{2-31-a}$), we
have
$$ d(Tx, Ty) \sim |\ta-t|+ |\xi-\eta|\sqrt{\sin (8\ta)\sin
(8t)}\sim d(x, y).$$ Thus, it follows by ($\ref{1-3}$) that
\begin{align*}
\rho_\al (Tx, Ty) &\sim \f {d(Tx, Ty)}\al +\f {
|\sqrt{\al-8\ta}-\sqrt{\al-8t}|} {\sqrt {\al}} \\
&\sim \f {d(x, y)}\al +\f { |\sqrt{\f \al8-\ta}-\sqrt{\f\al8-t}|}
{\sqrt {\al}}\sim \rho_{\al/8} (x, y),\end{align*} which proves
Assertion (i).

(ii)\   \ Assertion (ii) follows directly from Assertion (i).

(iii)\  \  Let $E$ be  a measurable subset of $ B(e,\al)$.  Then
using ($\ref{4-3}$), we obtain
$$ |E|=8\int_{T^{-1}(E)} D(x\cdot e)\, d\sa(x).$$
Assertion (iii) then follows by noticing  that $D(\cos\ta)\sim 1$
whenever $\ta \in [0, \f { \pi-\va}8]$.

(iv)\   \  Assertion (iv) is a simple consequence of Assertions
(ii) and (iii),  Lemma 2.2 (iii) and the fact that $\tr_{r, B(e,
\f \al 8)}(T^{-1}x)\sim\tr_{r, B(e,\al)}(x)$ for any $x\in
B(e,\al)$ and $r\in (0,1)$.\end{proof}

Now we are in a position to prove Theorem 1.1 and Corollaries
1.2--1.4
 in the case when $\al\in [\f 12, \pi-\va]$.

\vspace{5mm}

 {\it Proof of Theorem 1.1.}\ \  Suppose that $\Ld$ is $(\rho_{\al}, \f \da
n)$-separable in $B(e,\al)$.  It then follows by Lemma 4.2 (i)
that $T^{-1} (\Ld)$ is $(\rho_{_{\al/8}}, \f \da { n
C_5})$-separable in $B(e,  \al/8)$. Thus, for any $\be \ge 1$,
\begin{align*}
\sum_{\og\in\Ld} &\(\max_{ x, y\in B_{\rho_\al} (\og, \f \da n)}
|f(x)-f(y)|^p\) |B_{\rho_\al} (\og, \f \da n)|\\
&\leq C_5 \sum_{\og\in \Ld} \(\max_{u, v \in
T^{-1}\bl(B_{_{\rho_{\al}}} (\og, \f {\da} n)\br)}
|f(Tu)-f(Tv)|^p\)\Bl|T^{-1}\bl(B_{_{\rho_{\al}}}
(\og, \f \da n)\br)\Br|\\
&\leq C_5 \sum_{z\in T^{-1}(\Ld)}\( \max_{u, v \in
B_{\rho_{\al/8}} (z, \f {C_5\da} n)} |f(Tu)-f(Tv)|^p\)
|B_{\rho_{\al/8}} (z,
\f {C_5\da} n)|\\
&\leq (C \da)^p \int_{B(e, \al/8)} |f(Tx)|^p \, d\sa(x)\leq (C
\da)^p\int_{B(e, \al)} |f(x)|^p \, d\sa(x),\end{align*}where in
the first inequality we have used Lemma 4.2 (iii), in the second
inequality we have used Lemma 4.2 (ii), in the third inequality we
have used  the already proven case of Theorem 1.1 applied to
$B(e,\al/8)$ and the polynomial  $f(Tx)\in\Pi_{8n}^d$, and in the
last inequality we have used ($\ref{4-3}$).  This proves
($\ref{1-6}$).

\vspace{5mm}
 {\it Proof of Corollary 1.2.}\   \ Suppose $\Ld$ is a
maximal $(\rho_\al, \f \da n)$-separable  subset of $B(e,\al)$.
Then by Lemma 4.2 (i)--(ii), $T^{-1} (\Ld)$ is $(\rho_{\al/8}, \f
{C_5^{-1}\da} n)$-separable in $B(e,\al/8)$ and
$$ \bigcup_{\og\in \Ld} B_{_{\rho_{\al/8}} } (T^{-1}\og, \f{C_5\da
}n)=B(e,\al/8).$$ Thus, slightly modifying the proof of  Corollary
1.2 in  the case  $\al\in (0,\f 12]$ given in Section 3, we
conclude that  there exists a constant $\da_1\in (0, 1)$ depending
only on $d$ such that if $\da \in (0, C_5^{-1}\da_1)$ then there
exists a sequence of positive numbers $\mu_\og$, \   \  $\og
\in\Ld$ such that
$$ \mu_\og \sim \tr_{_{\f\da n, B(e, \al/8)}}
(T^{-1}\og)\sim \tr_{_{\f\da n, B(e,\al)}}(\og),\   \   \   \og
\in\Ld$$ and such that for any $P\in \Pi_{8(n+d)}^d$,
$$\int_{B(e, \al/8)} P(y)\, d\sa(y) =\sum_{\og\in\Ld} \mu_\og
P(T^{-1}\og).$$ It then follows by Lemma 4.1 that for any $f\in
\Pi_n^d$,
$$\int_{B(e, \al)} f(y)\, d\sa(y)=
8 \int_{B(e, \al/8)} f(Ty) D(y\cdot e)\, d\sa(y)=\sum_{\og\in \Ld}
\(8\mu_\og D(e\cdot T^{-1}\og )\) f(\og).$$ Now  setting
$$\ld_\og
=8\mu_\og D(e\cdot T^{-1}\og), \   \   \og\in\Ld$$ and noticing
that $D(x\cdot e)\sim 1$ for $x\in B(e, \al/8)$, we deduce
Corollary 1.2 with $\da_0=C_5^{-1} \da_1$.\hb

\vspace{5mm}

 {\it Proofs of Corollaries  1.3 and 1.4.}\   \  First,  note that given  $\be
\ge 1$
 and an arbitrary $(\rho_\al, \f\da n)$-separable  subset $\mathcal{A}$ of  $B(e,\al)$, we have,  for
any $x\in B(e,\al)$,
\begin{align} \sum_{\xi\in\mathcal{A}} \chi_{_{B_{_{\rho_\al}}(\xi,\be \da/n)}}
(x)& =\sum_{\xi\in\mathcal{A}}
\chi_{_{T^{-1}\bl(B_{_{\rho_\al}}(\xi,\be
\da/n)\br)}} (T^{-1}x)\notag\\
&\leq \sum_{\eta\in T^{-1}(\Ld)}\chi_{_{B_{_{\rho_{\al/8}}}(\eta,
C_5\be \da/n)}} (T^{-1}x)\leq C_\be,\label{4-4}\end{align} where
the first inequality follows by Lemma 4.2 (ii), and  the second
inequality follows by Lemma 2.2 (v) and Lemma 4.2 (i). Now the
rest of the proofs are almost identical to those for the case
$\al\in (0,\f 12]$. We omit the details. \hb

\section{Concluding remarks}
\subsection{ Weighted  inequalities on spherical caps}
Let $\al\in (0, \f12]$ and let $e$ be a fixed point on $\sph$. A
weight function $W$ on $B(e,\al)$ is called a doubling weight if
there exists a constant $L>0$, called doubling constant, such that
for every $x\in B(e,\al)$ and $r\in (0,1)$,
$$\int_{B_{\rho} (x, 2r)} W(y)\, d\sa(y)\leq L \int_{B_{\rho} (x, r)} W(y)\,
d\sa(y),$$ where $\rho\equiv \rho_{_{B(e,\al)}}$  is  defined by
($\ref{1-3}$). Associated with a weight function $W$ on
$B(e,\al)$, we define
$$ W_n (x) :=\f 1{ | B_{\rho} (x, \f 1n)|} \int_{B_{\rho} (x,
\f1n)} W(y)\, d\sa(y),\    \   n=1,2,\cdots,\    \   x\in
B(e,\al).$$ It   follows by ($\ref{2-29-a}$) that  for a doubling
weight $W$ on $B(e,\al)$,
\begin{equation}\label{5-1}
W_n(x)\leq C L( 1+ n \rho (x, y))^{2+\f{\ln L}{\ln 2}} W_n (y),\ \
\ \text{for any $x, y\in B(e,\al)$},\end{equation} where $C$ is a
constant  depending  only on $d$.

We have the following theorem:

\begin{thm} Let $e\in\sph$, $1\leq p<\infty$ and  $\al\in (0,\f12]$.  Let
$W$ be a doubling weight on $B(e,\al)$. Then for any
$f\in\Pi_n^d$,
\begin{equation}\label{5-2}
\int_{B(e,\al)} |f(x)|^p W(x)\, d\sa(x) \sim  \int_{B(e,\al)}
|f(x)|^p W_n(x)\, d\sa(x),\end{equation} where the constant of
equivalence depends only on $d$, $p$ and the doubling constant of
$W$. Moreover, there exists a constant $\da_0$ depending only on
$d$, $p$,  and the doubling constant of $W$ such that  for  any
maximal $(\f \da n,\rho)$-separable subset $\Ld$ of $B(e,\al)$
with $\da \in (0,\da_0)$,  and any $f\in\Pi_n^d$, we have
\begin{align}
\int_{B(e,\al)} |f(x)|^pW(x)\, d\sa(x) &\sim \sum_{\og\in\Ld}
\(\max_{x\in B_{\rho} (\og,\f \da n)} |f(x)|^p\) \int_{B_\rho
(\og, \da/n)}
W(y)\, d\sa(y)\label{5-3}\\
&\sim \sum_{\og\in\Ld} \(\min_{x\in B_{\rho} (\og,\f \da n)}
|f(x)|^p\) \int_{B_\rho(\og, \da/n)} W(y)\,
d\sa(y),\label{5-4}\end{align} where the constants of equivalence
depend only on $d$, $p$ and the doubling constant of $W$.
\end{thm}

\begin{proof} For simplicity, associated with  a function $f$ on $B(e,\al)$,
we define
$$\osc(f)(x, r) =\max_{y, z\in B_\rho (x, r)} |f(y)-f(z)|,\    \
x\in B(e,\al),\   \  r>0.$$

For the proof of Theorem 5.1, we   claim that it is sufficient to
prove that for any $(\f \da n, \rho)$-separable subset $\Ld$ of
$B(e,\al)$ and any $f\in\Pi_n^d$,
\begin{equation}\label{5-5}
\sum_{\og\in\Ld} |\osc (f)(\og,  \da/ n)|^p \int_{B_\rho(
\og,\f\da n)} W_n (y) \, d\sa(y)\leq (C_6 \da)^p \int_{B(e,\al)}
|f(x)|^p W_n(x)\, d\sa(x),\end{equation} where $C_6$ depends only
on $d$, $p$ and the doubling constant of $W$. In fact, once
($\ref{5-5}$) is proved,  then setting $\da_0= \f 1{4C_6}$, and
taking into account Lemma 2.2 (v), we conclude that  for  any
maximal $(\f {\da_0} n,\rho)$-separable subset $\Ld$ of $B(e,\al)$
and any $f\in\Pi_n^d$, we have
\begin{align}\int_{B(e,\al)} |f(x)|^pW_n(x)\, d\sa(x) &\sim
\sum_{\og\in\Ld}
 \(\max_{x\in B_{\rho} (\og, \f {\da_0} n)} |f(x)|^p\) \int_{B_\rho(\og, {\da_0}/n)}
W_n(y)\, d\sa(y)\label{5-6}\\
&\sim \sum_{\og\in\Ld}
 \(\min_{x\in B_{\rho} (\og, \f {\da_0} n)} |f(x)|^p\) \int_{B_\rho(\og, {\da_0}/n)}
W_n(y)\, d\sa(y).\label{5-7}
\end{align}
Equation ($\ref{5-2}$) then follows by ($\ref{5-6}$),
($\ref{5-7}$), Lemma 2.2 (v) and the doubling property of $W$. On
the other hand, if $\Ld$ is an arbitrary maximal $(\f {\da}
n,\rho)$-separable subset  of $B(e,\al)$ with $\da\in (0,\da_0)$,
then setting $n_1= n\da_0/\da$, applying ($\ref{5-2}$),
($\ref{5-6}$) and ($\ref{5-7}$) to $f\in \Pi_{n_1}^d$, and in view
of Lemma 2.2 (v) and  the doubling property of $W$, we deduce
Equations ($\ref{5-3}$) and ($\ref{5-4}$).

 Thus, it remains to prove ($\ref{5-5}$).  We sketch the proof  as follows.
 First, we note that by ($\ref{5-1}$) and the standard technique in
 [D1], there exists a sequence of
positive polynomials $Q_n\in \Pi_n^d$ on $B(e,\al)$  such that $
W_n\sim Q_n^p$ and \begin{equation*}\osc (Q_n) (x, \da/ n) \leq C
\da Q_n (x),\ \ \ \ x\in B(e,\al).\end{equation*} It then follows
that
 \begin{align*}
 W_n (\og) \(\osc(f)
(\og, \f \da n)\)^p \leq& C   \(\osc(fQ_n) (\og, \f \da n)\)^p+ C
\(\max_{y\in B_{\rho} (\og, \f \da n)}|f(y)|^p\)
 \(\osc(Q_n) (\og, \f \da n)\)^p\\
 \leq & C  \(\osc(fQ_n) (\og, \f \da n)\)^p+C \da^p  \max_{y\in B_{\rho} (\og, \f \da n)}
 |f(y)Q_n(y)|^p, \end{align*}which combined with Theorem 1.1 and Corollary 1.3
 implies the desired inequality
 ($\ref{5-5}$).
 This completes the proof.
\end{proof}

Finally, we conjecture  that ($\ref{5-5}$) with $W_n$ replaced by
$W$ remains true. Note that by Lemma 3.2, this conjecture is true
when $d=1$.

\subsection {Analogous results on  spherical collars.}

Let  $e\in\sph$ and $0<\al <\be \leq \pi$. Recall that
$$B(e;\al,\be) =\{ x\in \sph:\    \   \al \leq d(x, e) \leq
\be\}$$ denotes the spherical collar centered at $e$ of spherical
height $\be-\al$.   We assume that $0<\al<\be <\pi-\va$ and $
\al\sim \be -\al$, where $\va\in (0,1)$ is a   given absolute
constant. We shall  keep this assumption for the rest of this
subsection. Without this assumption, some of the statements below
may not be true.

Associated with the spherical collar $B(e; \al,\be)$,  we define
\begin{equation}\label{5-8-a}\rho_{_{B(e; \al,\be)}}(x, y): =\f 1\al \sqrt{ |x-y|^2 + \al |
\sqrt{b_x}-\sqrt{b_y}|^2},\     \    x, y\in B(e;
\al,\be),\end{equation} where $b_x\equiv b_{x, B(e;\al,\be)}$
denotes the shortest distance from $x\in B(e;\al,\be)$ to the
boundary of $B(e;\al,\be)$, that is
$$b_x\equiv b_{x, B(e;\al,\be)}:=\min\Bl\{ d(x, y):\   \    y\in\sph,
 \  d(y, e) =\al\   \  \text{or}\   \  d(y, e)=\be\Br\}.$$ It
is easily seen that $\rho_{_{B(e; \al,\be)}}$ is a metric on $B(e;
\al, \be)$.

For $x= \xi\sin \ta + e \cos \ta$ and $y=\eta\sin t+e \cos t$ with
$\xi,\eta \in \mathbb{S}_e^{d-1}$ and $\ta, t\in [\al,\be]$, we
define
$$ \rho_6(x, y) := \max \Bl\{ |\xi-\eta|,\    \   \rho_{ _{[\al, \be]} }
(\ta, t)\Br\},$$ where
$$
 \rho_{[\al,\be]} (\ta, t): = \f
1\al \sqrt{ |\ta-t| ^2 +\al |\sqrt{b_{_{\ta, [\al,\be]}}}
-\sqrt{b_{_{t, [\al,\be]}}}|^2}$$ and $  b_{_{u, [\al,\be]}}$
denotes the shortest distance from $u\in[\al,\be]$ to the boundary
of the interval $[\al,\be]$, that is, $$  b_{_{u, [\al,\be]}}: =
\min \Bl\{ | u-\al|,\   \ |u-\be|\Br\}.$$ It turns out that in the
case $\al<\f\pi2$, $\rho_{_{B(e;\al,\be)}}$ and $\rho_6$ are
equivalent on the whole spherical collar $B(e; \al,\be)$. (The
proof of this fact is similar to that of Lemma 2.2 (i).)

Now our main results can be stated as follows:

\begin{thm}
Let $ \da\in (0,1)$ and   $1\leq p<\infty$.  Let $\rho\equiv
\rho_{_{B(e;\al,\be)}}$ be defined by ($\ref{5-8-a}$) and let $\Ld$
be a $(\f \da n, \rho)$-separable subset of $B(e; \al,\be)$. Then
for all $f\in\Pi_n^d$, we have
$$ \sum_{\og\in\Ld} \(\max_{x, y\in B_{\rho} (\og, \f \da
n)}|f(x)-f(y)|^p\) \bl|B_{\rho} (\og,  \da/ n)\br| \leq (C \da)^p
\int_{B(e; \al, \be)} |f(x)|^p \, d\sa(x),$$ where  the constant $C$
depends only on $d$ and $p$.\end{thm}

\begin{cor} There exists a constant $\da_0\in (0,1)$ depending only on
$d$ such that for any $\da\in (0,\da_0)$ and any maximal $(\f \da
n, \rho_{_{B(e;\al,\be)}})$-separable subset $\Ld$ of
$B(e;\al,\be)$ there exists a sequence of positive numbers
$\ld_\og$, $\og\in\Ld$ for which the following cubature formula
holds for all $f\in\Pi_n^d$:
$$\int_{B(e;\al,\be)} f(y)\, d\sa(y)=\sum_{\og\in\Ld} \ld_\og
f(\og).$$\end{cor}

Results  similar to Corollaries 1.3 and 1.4 can also be deduced
from Theorem 5.2.

For the  proofs of Theorems 5.2 and Corollary 5.3, the equivalence
between the metrics $\rho_{_{B(e;\al,\be)}}$ and $\rho_6$ plays an
important role. Since the proofs run along the same lines as those
of Theorem 1.1 and Corollary 1.2 given in Section 3, we omit the
details.




\end{document}